\documentclass[10pt]{amsart}
\allowdisplaybreaks[4]
\usepackage{amssymb,color}
\usepackage[square,compress,comma, numbers,sort]{natbib}

\usepackage[colorlinks=true, citecolor=blue, linkcolor=blue]{hyperref}
\usepackage{amsthm}

\definecolor{c20}{rgb}{0.,0.7,0.}
\definecolor{c30}{rgb}{0.,0.,1.}
\definecolor{c40}{rgb}{1,0.1,0.7}
\definecolor{c50}{rgb}{1,0,0}
\definecolor{c60}{rgb}{0,0.9,0.1}

\newcommand{\kb}[1]{\boldsymbol{#1}}
\newcommand{\vk}[1]{\kb{#1}}

\newcommand{\ve}{\varepsilon}

\newcommand{\abs}[1]{\lvert #1 \rvert}
\newcommand{\ABs}[1]{ \biggl \lvert #1 \biggr \rvert}
\newcommand{\norm}[1]{\lVert #1 \rVert}

\newcommand{\EE}[1]{\mathbb{E}\left\{  #1 \right\}}

\newcommand{\pk}[1]{\mathbb{P} \left(#1 \right) }

\newcommand{\R}{\mathbb{R}}

\newcommand{\inr}{\in \R}

\newcommand{\ldot}{,\ldots,}

\newcommand{\limit}[1]{\lim_{#1 \to   \infty}}

\newcommand{\BQN}{\begin{eqnarray}}
\newcommand{\EQN}{\end{eqnarray}}
\newcommand{\BQNY}{\begin{eqnarray*}}
\newcommand{\EQNY}{\end{eqnarray*}}

\def\eL#1{\textcolor{c30}{#1}}
\def\eL#1{{#1}}
\def\eE#1{\textcolor{c20}{#1}}
\def\CE#1{\textcolor{c30}{#1}}
\def\CE#1{#1}
\def\eE#1{#1}
\def\K1#1{\textcolor{cyan}{#1}}
\def\K1#1{#1}
\def\kk#1{\textcolor{black}{#1}}

\newcommand{\BS}{\begin{sat}}
\newcommand{\ES}{\end{sat}}
\newcommand{\BT}{\begin{theo}}
\newcommand{\ET}{\end{theo}}
\newcommand{\BK}{\begin{korr}}
\newcommand{\EK}{\end{korr}}

\newcommand{\BD}{\begin{de}}
\newcommand{\ED}{\end{de}}

\newcommand{\BIT}{\begin{itemize}}
\newcommand{\EIT}{\end{itemize}}
\newcommand{\BDI}{\begin{description}}
\newcommand{\EDI}{\end{description}}

\newcommand{\BRM}{\begin{remarks}}
\newcommand{\ERM}{\end{remarks}}

\newcommand{\BEL}{\begin{lem}}
\newcommand{\EEL}{\end{lem}}

\newtheorem{theo}{Theorem}[section]
\newtheorem{sat}[theo]{Proposition}
\newtheorem{de}[theo]{Definition}
\newtheorem{lem}[theo]{Lemma}

\newtheorem{korr}[theo]{Corollary}
\newtheorem{remark}[theo]{Remark}
\newtheorem{remarks}[theo]{Remarks}

\newcommand{\nelem}[1]{{Lemma \ref{#1}}}

\newcommand{\netheo}[1]{{Theorem \ref{#1}}}

\newcommand{\prooftheo}[1]{ \textbf{Proof of Theorem} \ref{#1} }

\newcommand{\prooflem}[1]{\textbf{Proof of Lemma} \ref{#1}}

\newcommand{\COM}[1]{}

\newcommand{\QED}{\hfill $\Box$ \\}

\def\Ga{\gamma}

\def\rw{\rightarrow}
\def\RW{\Rightarrow}

\def\IF{\infty}

\def\LT{\left}
\def\RT{\right}

\def\H{\mathcal{H}}

\def\eG#1{\textcolor{black}{#1}}
\def\eHH#1{\textcolor{black}{#1}}
\def\eHHH#1{\textcolor{black}{#1}}
\def\eHH#1{#1}
\def\eHHH#1{#1}

\def\ehb#1{\textcolor{c50}{#1}}
\def\ehb#1{#1}
\def\ehc#1{\textcolor{c50}{#1}}
\def\ehc#1{#1}
\def\ehd#1{\textcolor{c50}{#1}}
\def\ehd#1{#1}

\def\KD#1{#1}
\def\ehhe#1{#1}

\def\EH#1{\textcolor{c40}{#1}}
\def\EH#1{#1}
\def\eeh#1{\textcolor{c40}{#1}}
\def\kd#1{\textcolor{cyan}{#1}}
\def\rd#1{\textcolor{red}{#1}}
\def\he#1{\textcolor{c50}{#1}}
\def\eeh#1{#1}
\def\kd#1{#1}
\def\rd#1{#1}
\def\he#1{#1}
\def\ehg#1{\textcolor{black}{#1}}
\def\blu#1{\textcolor{black}{#1}}
\def\blu#1{\textcolor{black}{#1}}

\def\vf{\sigma^2}
\def\barg{\overline{\gamma}}

\def\vfn{ V_\varphi}

\begin{document}

\title[ $\gamma$-reflected Gaussian processes with stationary increments]
{Extremes of $\gamma$-reflected Gaussian processes with stationary increments}

\author{Krzysztof D\c{e}bicki}
\address{Krzysztof D\c{e}bicki, Mathematical Institute, University of Wroc\l aw, pl. Grunwaldzki 2/4, 50-384 Wroc\l aw, Poland}
\email{Krzysztof.Debicki@math.uni.wroc.pl}
\author{Enkelejd Hashorva}
\address{Enkelejd Hashorva, Department of Actuarial Science, University of Lausanne, UNIL-Dorigny 1015 Lausanne, Switzerland}
\email{enkelejd.hashorva@unil.ch}
\author{Peng Liu}
\address{Peng Liu,
 Department of Actuarial Science, University of Lausanne, UNIL-Dorigny 1015 Lausanne, Switzerland and Mathematical Institute, University of Wroc\l aw, pl. Grunwaldzki 2/4, 50-384 Wroc\l aw, Poland}
\email{peng.liu@unil.ch}

\bigskip

\date{\today}
 \maketitle

\bigskip
{\bf Abstract:}
For a given centered Gaussian process with stationary increments $X(t), t\geq 0$ and $c>0$, let
$$
W_\gamma(t)=X(t)-ct-\gamma\inf_{0\leq s\leq t}\left(X(s)-cs\right), \quad t\geq 0$$
denote the $\gamma$-reflected process, where $\gamma\in (0,1)$. This process is
important for both queueing  and   risk theory. In this contribution we are concerned
with the asymptotics, as $u \to \IF$, of
$$\pk{\sup_{0\leq t\le  T}W_\gamma(t)> u }, \quad T\in (0,\IF].$$
Moreover, we investigate the approximations of first and last passage times for given large threshold $u$. We apply our findings to the cases with $X$ being
the multiplex fractional Brownian motion and the Gaussian integrated  process.
As a by-product we derive an extension of Piterbarg inequality for threshold-dependent random fields.\\



\section{Introduction}
The seminal contribution \cite{HP99} derived the exact asymptotics, as the initial
capital $u$ tends to infinity, of the ruin probability
$$ \CE{\psi_{0,\IF}}(u)=\pk{\sup_{t\ge 0} W_0(t) >u}, \quad W_0(t):= X(t)- ct, c>0$$
 for some general centered Gaussian processes $X(t), t\ge 0$.
 A key merit of the aforementioned paper is that it paved the way for the
  \eE{study} of the tail \eE{asymptotics} of supremum of Gaussian processes
   with trend over unbounded intervals.
With a strong impetus from \cite{HP99} a wide range of
asymptotic results for supremum of such threshold dependent
families of Gaussian processes were obtained in
\cite{DE2002,HP2004,DI2005,MR2275917,MR2441468,MR2827458,MR2738340,MR2462285}.

\kd{This paper is devoted to the analysis of extremes of $\gamma$-reflected processes $W_\gamma$, defined as
\BQNY
W_\gamma(t)=X(t)-ct-\gamma\inf_{0\leq s\leq t}\left(X(s)-cs\right), \quad \gamma\in[0,1),
\EQNY
where $X$ is a centered Gaussian process with stationary increments and $c>0$.
}
The analysis of $\gamma$-reflected processes is of interest for both queueing and risk theory.
In risk theory $\gamma$ is related to a fixed tax-payment rate, \kd{with}
\begin{eqnarray}
\psi_{\Ga,\IF}(u)=\CE{\pk{\inf_{0\leq t< \IF}\Bigl( u-W_\gamma(t)\Bigr) <0}}
\kd{=\pk{\sup_{0\leq t< \IF} W_\gamma(t)>u}}\label{psi.gamma}
\end{eqnarray}
\kd{representing the infinite-time ruin probability with initial capital $u$,}
see e.g.,  \cite{AA10}.
 For $\gamma=1$, $W_1$ has also interpretation as a transient queue length process
in a fluid queueing system
fed by $X$ and emptied with constant rate $c>0$,
see e.g., \cite{Harr85,AwaGly09,ZeeGly00,MandjesKrzys}.\\
More importantly, investigation of extremes of such processes is related to investigation of extremes of Gaussian random fields with interesting structures as already shown
in \cite{HA2013}. Therein the \kd{asymptotics of (\ref{psi.gamma})
for $X=B_H$ a fractional Brownian motion with Hurst index $H\in (0,1)$ has been investigated.}
Using the \kd{self-similarity of $B_H$,} for any $u>0$ and $X=B_H$ we have
\BQN
\psi_{\Ga,\IF}(u)
&=& \pk{\sup_{t\ge 0}\left(X(t)-ct-\gamma\inf_{s\in[0,t]}(X(s)-cs)  \right)>u}\nonumber\\
&=&\CE{\pk{\sup_{0\le s \le t<\IF} \frac{X(tu)- \gamma X(s u)  }{1+ c(t- \gamma s)} > u}}\notag \\
&=&\pk{\sup_{0\le s \le t<\IF} Y(s,t)  > u^{1- H}},
\label{Infruin}
\EQN
where $Y(s,t):=\frac{X(t)-\gamma X(s)}{1+ct-c\Ga s} $. In view of \eqref{Infruin} it is clear that for $X$ being an fBm, the approximation of $\psi_{\gamma,\IF}(u)$ \CE{as $u\to \IF$} is
 closely related to the study of supremum of the  Gaussian random field $Y$. The fact that $Y$ does not depend on the threshold $u$ is crucial \CE{and leads to substantial simplifications of the problem at hand.}
\CE{However, for a general centered Gaussian process $X$ with stationary increments, due to the lack of self-similarity, one has to analyse the tail behaviour of threshold-dependent random field
\BQN\label{maine}
Y_u(s,t)=
\frac{X(tu)-\gamma X(su)}{1+ct-c\Ga s}, \quad  s,t \in [0,\IF)
,\EQN
 which significantly increases the complexity of the  problem due to the explicit dependence
 on the threshold $u$.}
 \kd{We
 overcome this difficulty by deriving extensions of two classical
 results in extreme value theory of Gaussian processes. In particular,
Lemma \ref{PICKANDS} provides a uniform (with respect to local behavior of variance-covariance structure
of family of processes $X_u$) version of the celebrated
Pickands-Piterbarg lemma, as given in, e.g., Theorem D.2 in \cite{Pit96}.
Lemma \ref{lemP} extends Piterbarg inequality to
 threshold-dependent Gaussian random fields.
The generality of these findings makes them also applicable to other problems related with extremes of threshold-dependent families of Gaussian random fields.} \\
Under some conditions on the variance function $\vf$, assuming in particular that it is regularly varying  with index $2\alpha_0$ and $2 \alpha_\IF$ at 0 and $\IF$, respectively, \CE{our main result presented in \netheo{TH1} below} gives an asymptotic expansion of  $\psi_{\gamma,\IF}(u)$ as $u \to \IF$. \CE{It turns out that three different types of asymptotics of $\psi_{\gamma,\IF}(u)$ take place, mainly
determined by} the following limit  (which we assume to exist)
\BQN\label{limphi}
\varphi:= \limit{u} \frac{\vf(u)}{u} \in [0,\IF],
\EQN
where $\vf(t)=Var(X(t))$.
\K1{Interestingly, this trichotomy is tightly related with the
dependence structure of $X$.
For example, if $X= B_H$,
we can distinguish the case of $\varphi \in (0,\IF)$, i.e., $X$ is a standard Brownian motion,  $\varphi=0$ if $H \in (0,1/2)$ which is the well-known case of {\it short range dependent} fBm and $\varphi=\IF$ corresponding to $H \in(1/2,1]$, i.e., the case of {\it long range dependent} fBm}.\\
 Comparing our findings with those obtained for $\gamma=0$ in \cite{DI2005}, using $\sim$ to denote the asymptotic equivalence, we obtain the following {\it asymptotic tax equivalence} (derived for $X=B_H$ in \cite{HA2013})
\BQN\label{gammazero}
 \psi_{\gamma,\IF}(u)& \CE{\sim} &
\mathcal{P}_{\vfn }^{\barg} \psi_{0,\IF}(u) , \quad \barg:= (1- \gamma)/ \gamma, \quad \gamma \in (0,1)
\EQN
as $u \to \IF$, with
\BQN \label{varphi}
 \vfn =\frac{\sqrt{2}c\gamma}{\varphi}X, \quad \text{if } \varphi \in (0,\IF), \quad
\vfn= B_{\alpha_{\varphi}}, \text{ if  } \varphi \in \{0, \IF\}.
\EQN
In our notation
\BQNY
\mathcal{P}_{Z}^a =  \EE{ \sup_{t\in [0,\IF)} e^{ \sqrt{2} Z(t)- (1+ a) Var( Z(t) ) }}, \quad a>0
\EQNY
\KD{denotes the generalised Piterbarg constant, where $Z$
is a centered Gaussian process with stationary increments and continuous sample paths.} Note in passing that by Theorem 1.1 in \cite{Genna88} a.s.\ continuity of $Z$ at each $t\in [0,S]$ is
 equivalent to the sample-continuous assumption above.
 Further, the constants \KD{$ \mathcal{P}_{B_H}^a,$ with $ B_H$} a standard fBm, are known only for
\BQN
 \mathcal{P}_{B_{1/2}}^a=1+\frac{1}{a}\ \ \text{and}\ \ \mathcal{P}_{B_1}^a=\frac{1}{2}\left(1+\sqrt{1+\frac{1}{a}}\right),
\label{eqpp}
\EQN
see e.g., \cite{Pit96, MR1993262,LongB}. For general \KD{$H\in(0,1)$, bounds for  $ \mathcal{P}_{B_H}^a$} are derived in \cite{debicki2011extremes,LongB}.
\\
The asymptotics in (\ref{gammazero}) shows that the generalised Piterbarg constant governs
the relation between the two ruin probabilities corresponding to the model with tax and without tax, i.e.,
it defines what we call {\it the asymptotic tax equivalence}.
\CE{However, in view of \cite{HJ13,PengLith} we know that for the case $X=B_H$, the tax rate $\gamma$ does not influence
the limiting distribution of the first and the last passage times. We investigate these \KD{problems in more general}
model\EH{s} for $X$.
Define therefore} the first and last passage times of $W_\gamma$ given that the ruin occurs by
\BQN\label{tau}
(\tau_1^*(u),\tau_2^*(u))\stackrel{d}{=} (\tau_1(u),\tau_2(u))\Bigl\lvert(\tau_1(u)<\IF),
\EQN
where $$\tau_1(u)=\inf\{t\geq 0, W_\gamma(t)>u\} \ \ \mbox{and} \ \ \tau_2(u)=\sup\{t\geq 0, W_\gamma(t)>u\},
$$
with the convention that $\inf\{\emptyset\}=\IF$ and $\sup\{\emptyset\}=0$. \CE{Here $\stackrel{d}{=}$ stands for equality of the distribution functions.}\\
\KD{ Complementary, in this contribution we address also finite-time horizon counterparts of the introduced above problems. Namely}
\BQN\label{eq:finit}
\psi_{\eL{\gamma}, T}(u):=\mathbb{P}\left(\sup_{0\leq t\leq T} W_\gamma(t)>u\right)
\EQN
for any finite $T>0$ is analysed, extending partial results on \eL{ $\psi_{0,T}$ \CE{given} in \cite{DebickiRol02}}. Moreover, we shall deal also with the approximation of the conditional first passage time
$$\tau_1(u)\bigl\lvert(\tau_1(u)<T)$$
\CE{as $u\to \IF$
\KD{(see \netheo{THF2}),}
which shows that the approximating random variable is exponentially distributed.\\
\KD{The family of Gaussian processes $X$ with stationary increments, considered in this contribution,
covers general classes such  as}
\begin{itemize}
	\item  [A)] Multiplex fBm model, i.e., }
\BQNY
X(t)=\sum_{i=1}^n B_{H_i}(t),  \quad 
t\geq 0,
\EQNY
\CE{with $B_{H_i}$'s being independent fBm's;
	\item [B)] Gaussian integrated process model, that is the case where
$\CE{X(t)}=\int_0^tY(s)ds, t\geq 0$ with $Y$ being a centered stationary  Gaussian process with \he{a.s.}\ continuous sample paths.}
\end{itemize}
\eL{Organization of the paper: In Section 2  we present some preliminaries, followed by the main
results  for the approximation of $\psi_{\gamma,T}(u), T\in (0,\IF]$,
the approximating joint distribution for conditional scaled first and last passage times for $T \in (0,\IF]$.
 Section 3 {\he{is dedicated to applications related to model A) and B)} mentioned above.}
\he{For reader's convenience, we postpone all the proofs to Section 4; whereas some very technical claims are presented in Appendix.}}

\def\dotvT{\dot \vf(T)}
\section{Main Results}
In the rest of this paper $X(t),t\ge 0$ is a centered Gaussian process with stationary increments,
\he{a.s.}\ continuous sample paths and variance function \KD{$\vf(t)$}.
An canonical example is $X=B_H, H \in (0,1]$ for which we have $\vf(t)= t^{2H}$.
{For a given centered Gaussian process $Z$ with \he{a.s.}\ continuous sample paths set}
\BQNY
\mathcal{H}_Z[0,S]=\EE{\sup_{t\in[0,S]}
e^{\sqrt{2}Z(t)-Var(Z(t))}} 
\EQNY
{and define (whenever the limit exits) the generalised Pickands constant ${\mathcal{H}_Z}$ by}
\BQNY
\mathcal{H}_Z=\lim_{S\rw\IF} S^{-1}\mathcal{H}_Z[0,S].
\EQNY
 See \cite{PicandsA,PicandsB, Pit96, MR1772400, HH2000, DE2002,DI2005,DE2014,DiekerY,Harper2, DEJ14,Tabis,Harper3, DM,mi:17,SBK,Shao, high1} for various definitions, existence and basic properties
  of Pickands constant.

\subsection{Infinite-time horizon}
\KD{First we focus on \he{the} infinite-time horizon case.}
Due to the stationarity of increments, the covariance of $X$ is directly
defined by $\vf$, therefore our assumptions on $X$ shall be reduced to assumptions on the variance function, namely:

{\bf AI}: $\vf(0)=0$ and $\vf(t)$ is regularly varying at $\IF$ with index $2\alpha_\IF\in(0,\eHH{2)}$.
Further, $\vf(t)$ is twice continuously differentiable on $(0,\IF)$ with its first derivative
\KD{$\dot{\vf}(t):=\frac{{\rm d} \sigma^2}{{\rm d}t}\left(t\right)$}
and second derivative
\KD{
$\ddot{\vf}(t):=\frac{{\rm d^2} \sigma^2}{{\rm d}t^2}\left(t\right)$
}
being ultimately monotone at $\IF$.\\
{\bf AII}: $\vf(t)$ is regularly varying at $0$ with index $2\alpha_0\in(0,2]$ and its first derivative $\dot{\vf}(t)$ is ultimately monotone as $t\rw 0$.\\
{\bf AIII}: $\vf(t)$ is increasing and $\frac{\vf(t)}{t^2}$ is decreasing over $(0,\IF)$.\\

\COM{
\BRM i) From {\bf AI-AII}, one can get that there exists a positive constant $C>0$ such that in a neighborhood of zero,
\BQNY
\vf(t)\geq C_1t^2
\EQNY
holds. (see Lemma 5.1 in \cite{}). Further, we know that
\BQN\label{A1}
g(t):=\frac{t^2}{\vf(t)},\ \ t\in(0,\IF)
\EQN
is a regularly varying function at infinity with index $2(1-\alpha_\IF)>0$ and bounded in a neighborhood of zero.\\
ii) It follows from {\bf AI} and {\bf AIII} that
\BQN\label{A2}
g_1(t):=\frac{t}{\dot{\vf}(t)},\ \ t\in(0,\IF)
\EQN
is a regularly varying function at infinity with index $2(1-\alpha_\IF)>0$ and bounded in a neighborhood of zero.
\ERM
}
\def\LU{\Lambda(u)}
\def\CU{\frac{ 1}{c} \sqrt{\frac{2\alpha_\IF \pi}{1-\alpha_\IF}}  }
Define $\varphi$ by \eqref{limphi} assuming that the limit exists. For notational simplicity we set
$$\eL{t_*}=\frac{\alpha_\IF}{c(1-\alpha_\IF)} \he{> 0}$$
and
\BQN\label{Delta}
\Delta_\gamma(u)=\left\{
\begin{array}{cc}
\overleftarrow{\sigma}
 \left(\frac{\sqrt{2}\sigma^2(u\eL{t_*})}{\gamma u(1+ct_{*})}\right),& if \ \varphi=\IF \ or\  0,\\
 1,&if \ \varphi\in(0,\IF),\\
 \end{array}
 \right.
\EQN
where $\overleftarrow{\sigma}$ is the {asymptotic} inverse of $\sigma$ {(see e.g., \cite{Res,Soulier} for details)}.

Let $t_u$ be \he{a} maximizer of \eL{$\frac{\sigma(ut)}{1+ct}$ over $t\ge 0$.}
\EH{In view of} Lemma \ref{L1} for $u$ large enough $t_u$ is unique and
$$\limit{u}t_u=t_* {\in (0,\IF)}.$$
\kd{Hereafter $\Psi$ stands for the survival function of an $N(0,1)$ random variable.
\\
Before stating our main result, let us
observe that
$$\psi_{\gamma,\IF}(u)=\pk{\sup_{0\leq s\leq t<\IF} \frac{X(tu)-\gamma X(su)}{u(1+ct-c\Ga s)}>1}$$
\he{is valid for any $u>0.$} Typically the most likely point to
reach high value $u$ for
a centered Gaussian random field corresponds to the point that maximizes
its variance function, i.e., in our case
$$(s_u, t_u):=\mbox{argsup}_{(s,t): \ 0\leq s\leq t<\IF} Var\left(\frac{X(tu)-\gamma X(su)}{u(1+ct-c\Ga s)}\right).$$
It will be  shown in Lemma \ref{L1}
that $s_u=0$ for $u$ large and thus
$t_u=\mbox{argsup}_{t:\ t\geq 0}\frac{\sigma(ut)}{u(1+ct)}$.
This  explains the exponential term in the derived asymptotics.\\
The following theorem extends results derived in \cite{HA2013}, where
the special case $X=B_H$ is considered.}

\BT\label{TH1} If {\bf AI-AIII} are satisfied, then for any $\gamma\in (0,1)$ and $\eL{\varphi \in [0,\IF]}$
we have
\BQNY
\psi_{\gamma,\IF}(u)\sim \CU \mathcal{H}_{\vfn} \mathcal{P}_{\vfn} ^{\barg } \frac{\sigma(u\eL{t_*})}{
\Delta_1(u)}\eG{\Psi\left(\frac{u(1+ct_u)}{\sigma(ut_u)}\right)},
\EQNY
with $\vfn= \frac{\sqrt{2}c}{\varphi}X $ if $\varphi\in (0,\IF)$ , $\vfn=B_{\alpha_\varphi}$ if $\varphi\in \{0,\IF\}$ and $\barg:= (1- \gamma)/\gamma.$
\ET
\ehb{An immediate application of the above theorem, together with the known results in \cite{DI2005} for the case
$\gamma=0$, yields
\KD{that, as $u\to \IF$}
$$ \psi_{\gamma,\IF}(u) \sim  
\mathcal{P}_{\CE{\vfn}}^{\barg} \psi_{0,\IF}(u).
$$
}
\ehc{The above asymptotic tax equivalence shows that $\psi_{\gamma,\IF}(u)$ is proportional to
$\psi_{0,\IF}(u)$ as $u\to \IF$, where the proportionality constant is determined by the generalised Piterbarg constant
$\mathcal{P}_{\CE{\vfn}}^{\barg} $.}

\BT\label{THA} If {\bf AI-AIII} are satisfied, then for any $\gamma\in (0,1)$ and $\eL{\varphi \in [0,\IF]}$
we have \ehb{the convergence in distribution}
\BQNY
\left(\frac{\tau_1^*-ut_u}{A(u)},\frac{\tau_2^*-ut_u}{A(u)}\right)\stackrel{d}{\rw} \left(\mathcal{N}, \mathcal{N}\right), \quad u\rw\IF,
\EQNY
  where
$A(u)=\frac{\sigma(u\eL{t_*})}{c}\sqrt{\frac{\alpha_\IF}{1-\alpha_\IF}}$ and \blu{$\mathcal{N}\sim N(0,1)$.}
\ET
\eG{The above result implies that the standardized conditional first pasage time $\frac{\tau_1^*-ut_u}{A(u)}$ and  last passage time $\frac{\tau_2^*-ut_u}{A(u)}$ both weakly converge to standard normal random variables and $\frac{\tau_2^*(u)-\tau_1^*(u)}{A(u)}\rw 0$ in probability as $u\to \IF$.}

\subsection{Finite-time horizon}
\eL{Next,  we consider the finite-time horizon ruin probability, investigating $\psi_{\gamma, T}$ for $T$ a finite positive constant.} \\
\eL{\eG{Since we consider the finite-time horizon}, we shall impose weaker assumptions on the variance function $\vf$, namely:}
\def\vfp{\eL{\dot{\vf}}}
\def\vfpp{\eL{\ddot{\vf}}}

{\bf BI}: $\sigma^2(0)=0$ and $\vf(t)$ is twice differentiable over interval $(0,T]$.\\
{\bf BII}: $\vf(t)$ is regularly varying at $0$ with index $2\alpha_0\in (0,2]$.\\
{\bf BIII}: \eG{For $t\in (0,T]$, the first derivative $\dot{\vf}(t)>0$  and $\frac{\vf(t)}{t^2}$ is decreasing.}

For notational simplicity we set below
$$q(u)=\overleftarrow{\sigma}\left(\frac{\sqrt{2}\vf(T)}{u+cT}\right).$$

\BT\label{THF1} Suppose that {\bf BI--BIII} hold and \eG{$\gamma \in (0,1)$}.\\
i) If $s=o(\vf(s))$ as $s\rw 0$,  \eL{then}
\BQNY
\psi_{\gamma,T}(u)&\sim &\mathcal{H}_{B_{\alpha_0}}\mathcal{P}_{B_{\alpha_0}}^{\barg}
\frac{2\sigma^4(T)}{\vfp(T)q(u)u^2}\Psi\left(\frac{u+cT}{\sigma(T)}\right).
\EQNY
ii) If $\lim_{s\to 0} \frac{\vf(s)}{s}=b\in (0,\IF)$, then
\BQNY
\psi_{\gamma,T}(u)&\sim& \mathcal{P}_{B_{1/2}}^{\frac{\vfp(T)}{b}}\mathcal{P}_{B_{1/2}}^{\beta(b,\gamma)}
\Psi\left(\frac{u+cT}{\sigma(T)}\right), \quad \beta(b,\gamma):=\frac{b(\gamma-\gamma^2)+\gamma\vfp(T)}{b\gamma^2}.
\EQNY
iii) If $\vf(s)=o(s)$ as $s\rw 0$, then
\BQNY
\psi_{\gamma,T}(u)&\sim& \Psi\left(\frac{u+cT}{\sigma(T)}\right).
\EQNY
\ET
\BRM\label{RE}
i) From the proof of Theorem \ref{THF1}, we can similarly get
\KD{the asymptotics of $\psi_{0,T}(u)$ (see also \cite{DebickiRol02}),
which compared with $\psi_{\gamma,T}(u), \gamma\in (0,1)$,  gives
}
$$
\psi_{\gamma,T}(u)\sim \mathcal{K}\psi_{0,T}(u),\ \ u\rw\IF,
$$
with
\BQNY
\mathcal{K}=\left\{\begin{array}{cc}
\mathcal{P}_{B_{\alpha_0}}^{\barg},& if s=o(\vf(s)),\\
\mathcal{P}_{B_{1/2}}^{\beta(b,\gamma)  },& if \lim_{s\to 0} \frac{\vf(s)}{s}=b\in (0,\IF),\\
1,& if \vf(s)=o(s).
\end{array}\right.
\EQNY
ii)
\kk{
The approach used in the proofs of Theorem \ref{TH1} and Theorem \ref{THF1}
enables us to find exact asymptotics of
$\psi_{\gamma, T_{u}}(u)$ as $u\to\infty$, for some scenarios where
the time-horizon $T_u$ is a deterministic function of $u$.
For example, if
$ut_u=o(T_u)$ as $u\to\infty$, then by the proof of Theorem \ref{TH1}
we have
$\psi_{\gamma, T_{u}}(u)\sim\psi_{\gamma, \infty}(u)$, $u\to\infty$.
Additionally, if
$T_{u}\rw T$ as $u\rw\IF$, then the asymptotics of
 $\psi_{\gamma, T_{u}}(u)$ \eL{can be obtained by replacing   $T_{u}$ with  $T$ in}  the corresponding formulas of Theorem \ref{THF1}.
On the other side,
the case $T_u\sim t^*u$ as $u\to\infty$,
is out of the approach given in this paper. We suspect that it
leads to the asymptotics of
qualitatively other type than derived in Theorems \ref{TH1}, \ref{THF1}.
}
\ERM
\blu{Next we consider a finite-time counterpart of Theorem \ref{THA}. While for the infinite-time horizon the limit distribution in Theorem \ref{THA} is Gaussian, as shown below, this is not the case for finite-time horizon, where the limit distribution is exponential. The intuitive explanation for this is that
the local behaviour of variance function of the considered Gaussian field
in neighbourhood of the variance maximizer plays the key role for the type of
the limit distribution. In particular, if the first derivative of the
variance function is positive at this point, then the limiting distribution is exponential,
while the first derivative equal to $0$ at that  point leads to limit with Normal distribution; compare Lemma \ref{L1} with Lemma \ref{Var}.}

\BT\label{THF2}
If {\bf BI--BIII} are satisfied and $\lim_{s\rw 0}\frac{\vf(s)}{s}\in[0,\IF]$, then \eE{the convergence in distribution}
\BQNY
\frac{\vfp(T)}{2\sigma^4(T)} u^2(T-\tau_1)\bigl\lvert (\tau_1\leq T)\stackrel{d}{\rw} \mathcal{E}
\EQNY
\eE{holds, as $u\to \IF$,} \eL{with $\mathcal{E}$ a unit exponential random variable.}
\ET

\section{Applications}
\eHH{In this section, we shall focus on two important classes of processes with stationary increments. We discuss first} the sum of independent \eE{fBm's} with different Hurst parameters and then investigate Gaussian integrated processes.

 \subsection{Multiplex fBm} Let next $B_{H_i}, 1\leq i\leq n$ be  independent standard fBm's with index $0<H_1<H_2\leq \cdots \leq H_{n-1}< H_n<1$ and define  for $t\ge 0$

\BQN\label{HFBM}
\K1{
X(t)= M_{\vk{H}}(t):=\sum_{i=1}^nB_{H_i}(t),  \quad \vk{H}=(H_1, \cdots, H_n).
}
\EQN
A motivation to consider such a process stems from the insurance models with tax,
where $\eL{B}_{H_i}$ represents the aggregated claims of the sub-portfolios of the insurance company. {We have that}
$$\vf(t)= \sigma_{M_{\vk{H}}}^2(t)=\sum_{i=1}^n t^{2H_i}$$
 satisfies {\bf AI-AIII} with $\alpha_0=H_1$, $\alpha_\IF=H_n$. Further,
\BQNY
\varphi=\left\{
\begin{array}{cc}
 \IF,& 1/2<H_n<1,\\
 1,&H_n=1/2,\\
 0,&0<H_n<1/2
 \end{array}
 \right.
\EQNY
{implying} the following result:

\BK
\K1{Suppose that $X$ is defined by
(\ref{HFBM}).}\\
i) If $0<H_n<1/2$, then
\BQNY
\psi_{\gamma,\IF}(u)
&\sim&\mathcal{H}_{B_{H_1}}\mathcal{P}_{B_{H_1}}^{ \barg }2^{\frac{H_1-1}{2H_1}}\sqrt{\pi}c^{\frac{2H_n-H_1H_n-H_1}{H_1}}H_n^{\frac{H_1-4H_n+2H_1H_n}{2H_1}}(1-H_n)^{\frac{4H_n-H_1-2H_1H_n-2}{2H_1}}\\
&& \times u^{\frac{H_1H_n-2H_n+1}{H_1}}\Psi\left(\inf_{t> 0}\frac{u(1+ct)}{\sigma_{M_{\vk{H}}}(ut)}\right).
\EQNY
ii) If $H_n=1/2$, then
\BQNY
\psi_{\gamma,\IF}(u)\sim \mathcal{H}_{\sqrt{2}cM_{\vk{H}}}\mathcal{P}_{\sqrt{2}c \gamma M_{\vk{H}}}^{ \barg }\sqrt{\frac{2\pi u}{c^3}}\Psi\left(\inf_{t> 0}\frac{u(1+ct)}{\sigma_{M_{\vk{H}}}(ut)}\right).
\EQNY
iii) If $1/2<H_n<1$, then
\BQNY
\psi_{\gamma,\IF}(u)
&\sim& \mathcal{H}_{B_{H_n}}\mathcal{P}_{B_{H_n}}^{ \barg }2^{\frac{H_n-1}{2H_n}}\sqrt{\pi}c^{1-H_n}H_n^{\frac{2H_n-3}{2}}(1-H_n)^{\frac{3H_n-2-2H_n^2}{2H_n}}u^{\frac{(1-H_n)^2}{H_n}}\Psi\left(\inf_{t> 0}\frac{u(1+ct)}{\sigma_{M_{\vk{H}}}(ut)}\right).
\EQNY
\EK
Moreover, since {\bf BI-BIII} are satisfied for $M_{\vk{H}}(t)$, we obtain for any $T>0$.
\BK
\K1{Suppose that $X$ is defined by (\ref{HFBM}).}\\
i) If $0<H_1<1/2$, then
\BQNY
\psi_{\gamma,T}(u)\sim \mathcal{H}_{B_{H_1}}\mathcal{P}_{B_{H_1}}^{\barg }2^{-\frac{1}{2H_1}}\frac{\left(\sum_{i=1}^nT^{2H_i}\right)^{\frac{2H_1-1}{H_1}}}
{\sum_{i=1}^nH_iT^{2H_i-1}}u^{\frac{1-2H_1}{H_1}}\Psi\left(\frac{u+cT}{\sqrt{\sum_{i=1}^nT^{2H_i}}}\right).
\EQNY
ii) If $H_1=1/2$, then
\BQNY
\psi_{\gamma,T}(u)\sim \mathcal{P}_{B_{1/2}}^{2\sum_{i=1}^nH_iT^{2H_i-1}}\mathcal{P}_{B_{1/2}}^{\frac{\gamma-\gamma^2+2\gamma\sum_{i=1}^nH_iT^{2H_i-1}}{\gamma^2}}
\Psi\left(\frac{u+cT}{\sqrt{\sum_{i=1}^nT^{2H_i}}}\right).
\EQNY
iii) If $1/2<H_1<1$, then
\BQNY
\psi_{\gamma,T}(u)\sim
\Psi\left(\frac{u+cT}{\sqrt{\sum_{i=1}^nT^{2H_i}}}\right).
\EQNY
\EK
\begin{remark}
\kd{In the above corollaries, the main contribution to the asymptotics depends on all $H_i$'s while the polynomial terms depend on the properties of variance function at time $0$ and $\IF$ which is determined by  Hurst parameters $H_1$ and $H_n$.
It follows from the fact that the formula under $\Phi(\cdot)$
comes from global optimum of the variance function of the appropriate Gaussian field, while
the polynomial part of the asymptotics
follows from the asymptotic relation between local behavior of variance and correlation
in the neighbourhood of the variance optimizer.
}
\end{remark}

\subsection{ Gaussian integrated  processes}  Suppose that
\BQN\label{SUP}
\CE{X(t)}=\int_0^tY(s)ds, t\geq 0,
\EQN
 where $Y$ is a stationary
 \KD{centered}
 Gaussian process with unit variance and  \he{a.s.}\ continuous sample paths.
Let \eeh{$R(t)=Cov\left(Y(s), Y(s+t)\right), s, t\geq 0$}. In this subsection, we shall consider two scenarios:\\
 {\bf SRD} (short-range dependent), i.e., we shall assume that\\
 i) $R(t)$ is decreasing over $[0,\IF)$,\\
 ii) 
 \KD{$\int_0^\IF R(t)dt=G \in \ehhe{(0,\IF)}$}.\\
 {\bf LRD} (long-range dependent), {i.e., we shall suppose that} \\
 i)   $R(t)$ is decreasing over $[0,\IF)$,\\
 ii) $R(t)$ is regularly varying at infinity with index $2H-2$, $H\in (1/2, 1)$.\\
It follows that {\bf AI-AIII} are satisfied if  $X$ is {\bf SRD} or {\bf LRD},
implying our next results.
 \BK Suppose that $X$ is defined by (\ref{SUP}).\\
 i) If $X$ is {\bf SRD}, then
 \BQNY
\psi_{\gamma,\IF}(u)\sim \mathcal{H}_{\frac{\sqrt{2}c}{G}X}\mathcal{P}_{\frac{\sqrt{2}\gamma c}{G}X}^{ \barg }\sqrt{\frac{2\pi G u}{c^3}}\Psi\left(\inf_{t> 0}\frac{u(1+ct)}{\sigma(ut)}\right).
\EQNY
ii) If $X$ is {\bf LRD}, then
\BQNY
\psi_{\gamma,\IF}(u)
&\sim& \mathcal{H}_{B_{H}}\mathcal{P}_{B_H}^{ \barg }2^{\frac{H-1}{2H}}\sqrt{\pi}c^{1-H}H^{\frac{1-4H+2H^2}{2H}}(1-H)^{\frac{3H-2-2H^2}{2H}}(2H-1)^{\frac{1-H}{2H}}\frac{u\sqrt{R(u)}}{\mathcal{R}^*(u)}\Psi\left(\inf_{t> 0}\frac{u(1+ct)}{\sigma(ut)}\right),
\EQNY
with $\mathcal{R}^*$ the {asymptotic} inverse function of $u\sqrt{R(u)}$.
 \EK

Since, {\bf BI-BIII} are satisfied (note that $\sigma^2(t)\sim t^2 =o(t)$ as $t\rw 0$) for $R(t)$  positive and decreasing on $[0,T]$,  applying Theorem {\ref{THF1}} we arrive at
\kk{the following corollary.}

\BK If $X$ is defined by (\ref{SUP}) with $R(t)$ positive and  decreasing on $[0, T]$, then
\BQNY
\psi_{\gamma,T}(u) &\sim& \Psi\left(\frac{u+cT}{\sigma(T)}\right), \quad u\to \IF.
\EQNY
\EK

\section{Proofs}
We \kk{begin with introduction of some useful notation. Namely} we write
	$$D:=\{(s,t): 0\leq s\leq t< \IF\}, \quad
\sigma_\gamma^2(s,t):=Var( X(t)-\gamma X(s)),$$ $$\sigma_{\gamma,u}(s,t):=\frac{\sigma_\gamma(us,ut)}{1+c(t-\gamma s)}$$
and set further  for $(s,t), (s_1,t_1)\in D$
$$r_u(s,t,s_1,t_1):=\eHH{Cor(} X(ut)-\gamma X(us), X(ut_1)-\gamma X(us_1)).$$
 Hereafter, $Q$, $Q_i, i=1, 2, \dots$ are positive constants that
 may change from line to line.
 For any non-zero random variable $X$ \EH{we shall define}  $$\overline{X}:=\frac{X}{\sqrt{Var(X)}}.$$
\he{In our proofs multiple limits appear;  the order when passing to limit is important.
We shall write for instance
$$ a_u(S,S_1) \to 0, \quad u\to \IF, S\to \IF, S_1 \to \IF$$
 to mean that
 $$\lim_{S_1\to \IF} \lim_{S\rw\IF}\lim_{u\to \IF} a_u(S,S_1)=0.$$
 This convention applies for other instances of double or triple limits.}

\rd{We briefly comment on some useful  properties of $\sigma$.
\CE{For  $\CE{\lambda}\in \mathbb{R}$, by {\bf AI} and {\bf AII}, the \KD{function}}
\BQN\label{GGAMMA1}
g_{\CE{\lambda}}(t):=\frac{\vf(t)}{t^{\CE{\lambda}}}
\EQN
is regularly varying at 0  with index $2\alpha_0-\CE{\lambda}$ and at infinity with index $2\alpha_\IF-\CE{\lambda}$.\\
Further, by uniform convergence theorem (UCT) in \cite{BI1989,EKM97,Soulier}, we have that for any $T>0$ and $0<\lambda<\min (2\alpha_0,2\alpha_\IF)$
\BQNY
\lim_{u\rw 0}\sup_{t\in (0,T]}\left|\frac{g_\lambda(ut)}{g_\lambda(u)}-|t|^{2\alpha_0-\lambda}\right|=0
\EQNY
implying that for any $T>0$, when $u$ \blu{is} sufficiently small
\BQN\label{re0}
\frac{\sigma^2(ut)}{\sigma^2(u)}=\frac{g_\lambda(ut)}{g_\lambda(u)}|t|^\lambda\leq 2|T|^{2\alpha_0-\lambda}|t|^\lambda, \quad t\in[0,T].
\EQN
Moreover,  Potter's bounds (see e.g., \cite{BI1989, Soulier,EKM97}) show that for any $0<\epsilon<2\alpha_0$, there exists $T>0$ and $Q_1, Q_2>0$ such that for all $0<s,t<T$
\BQN\label{re00}
Q_1\min\left(\left(\frac{t}{s}\right)^{2\alpha_0-\epsilon},\left(\frac{t}{s}\right)^{2\alpha_0+\epsilon} \right)\leq\frac{\sigma^2(t)}{\sigma^2(s)}\leq Q_2\max\left(\left(\frac{t}{s}\right)^{2\alpha_0-\epsilon},\left(\frac{t}{s}\right)^{2\alpha_0+\epsilon} \right).
\EQN}
\subsection{\prooftheo{TH1}}
\he{First, we re-write for any $u>0$ the ruin probability $\psi_{\gamma,\IF}(u)$ as}
\BQNY
\psi_{\gamma,\IF}(u)&=& \pk{\sup_{t\ge 0}\left(X(t)-ct-\gamma\inf_{s\in[0,t]}(X(s)-cs)  \right)>u}\nonumber\\
&=&\CE{\pk{\sup_{0\le s \le t<\IF} \frac{X(tu)- \gamma X(s u)  }{1+ c(t- \gamma s)} > u}}\notag \\
&=&\pk{\sup_{(s,t)\in D} Z_u(s,t) > m(u)},
\EQNY
with
\BQN\label{ZU}
Z_u(s,t)=\left(\frac{X(ut)-\gamma X(us)}{1+c(t-\gamma s)}\right)\left(\frac{1+ct_u}{\sigma(ut_u)}\right), \quad \ (s,t)\in D, u>0, \quad 
\EQN
and
\BQN\label{MU}
m(u)=\inf_{t\geq 0}\frac{u(1+ct)}{\sigma(ut)}=\frac{u(1+ct_u)}{\sigma(ut_u)}, \quad u>0.
\EQN
Hereafter we shall denote
\BQN\label{EU}
\quad E(u):= E_1(u)\times E_2(u),\quad
E_1(u)=\left[0,\frac{\overleftarrow{\sigma}\left(u^{-1}\sigma^2(u)\ln u\right)}{u}\right), \quad E_2(u)=\left(t_u-\frac{\sigma(u)\ln u}{u}, t_u+\frac{\sigma(u)\ln u}{u}\right).
\EQN
 \K1{
As it will be shown below, \he{the} set $E(u)$ covers
sufficiently large neighbourhood of the maximizer of variance of $Z_u$
in order to determine
the asymptotics of $\psi_{\gamma, \IF}(u)$ by supremum of  $Z_u(s,t)$ over $E(u)$.}
More formally, for any $u>0$ we write
\BQN\label{P1}
\Theta(u)\leq \psi_{\gamma,\IF}(u)\leq \Theta(u)+
\Theta_0(u),
\EQN
with
\BQNY
\Theta(u)=\pk{\sup_{(s,t)\in E(u)}Z_u(s,t)>m(u)}, \quad \rd{\Theta_0(u)=\pk{\sup_{(s,t)\in D \eE{\setminus} \CE{E(u)}}Z_u(s,t)>m(u)}.}
\EQNY
\blu{The strategy of the proof is to
derive first  the exact asymptotics of $\Theta(u)$ as $u\to \infty$ and then to show that
 (recall (\ref{P1})) that $ \limit{u} \Theta_0(u)/\Theta(u)=0$.}

\kd{Before proceeding to details of these steps of the proof, we
summarize some dependence properties of the analyzed Gaussian field which will be needed in our proofs.}

\subsubsection{Dependence structure of $Z_u$}

\kd{Proofs of the following lemmas are deferred to Appendix.}
\BEL\label{L1}
\eHH{If the variance function $\vf$ of $X$ satisfies {\bf AI-AII},  then}
for $u$ large enough, the unique maximizer of $\sigma_{\gamma,u}(s,t)$ \CE{over $D$} is  $(0,t_u)$ and $\limit{u}
t_u=\eL{t_*} \CE{\in (0,\IF)}$. Moreover, for any $0<\epsilon<\min (a_1,a_2)$, when $u$ \blu{is} large enough and $\delta$ \blu{is} small enough
\BQNY
(a_1-\epsilon)(t-t_u)^2+(a_2-\epsilon)\frac{\vf (us)}{\vf (u)}\leq 1-\frac{\sigma_{\gamma,u}(s,t)}{\sigma_{\gamma,u}(0,t_u)}\leq (a_1+\epsilon)(t-t_u)^2+(a_2+\epsilon)\frac{\vf (us)}{\vf (u)}, \ \ |t-t_u|<\delta, 0\leq s<\delta,
\EQNY
with
$$  a_1=:\frac{c^2(1-\alpha_\IF)^3}{2\alpha_\IF}, \quad
 a_2=:\frac{\gamma(1-\gamma)}2 \Biggl[\frac{c(1-\alpha_\IF)}{ \alpha_\IF}\Biggr]^{2\alpha_{\IF}}.$$
\EEL
\BEL\label{L2}
If {\bf AI-AIII} are satisfied and \CE{$\delta_u>0, u>0$ are such that} $\limit{u}\delta_u= 0$, then we have
\BQNY
\lim_{u\rw \IF}\sup_{(s,t)\neq (s_1,t_1)\in [0,\delta_u)\times(t_u-\delta_u,t_u+\delta_u)}\left|\frac{1-r_u(s,t,s_1,t_1)}{\frac{\vf (u|t-t_1|)+\gamma^2\vf (u|s-s_1|)}{2\vf (u\eL{t_*})}}-1\right|=0.
\EQNY
\EEL

\rd{\subsubsection{\he{Asymptotic upper bound for} $\Theta_0(u)$}
For notational simplicity we define next (recall that \rd{$D=\{(s,t): 0\leq s\leq t<\IF\}$})
 $$D_T=\{(s,t): 0\leq s\leq t\leq T\}, \quad D_T^c=D\setminus D_T,  \quad D_{\delta,u}=D_T \ehc{\setminus}([0,\delta]\times[t_u-\delta,t_u+\delta])$$
   and
   $$D_{\delta,u}^*=([0,\delta]\times[t_u-\delta,t_u+\delta]) {\setminus}E(u).$$
For any $u>0$
\BQNY\label{NEG1}
&&\pk{\sup_{(s,t)\in \eE{D \setminus } E(u)}Z_u(s,t)>m(u)}\nonumber\\
&\quad&\quad\leq \pk{\sup_{(s,t)\in D_T^c}Z_u(s,t)>m(u)}+\pk{\sup_{(s,t)\in D_{\delta,u}}Z_u(s,t)>m(u)}+\pk{\sup_{(s,t)\in D_{\delta,u}^*}Z_u(s,t)>m(u)}\nonumber\\
&\quad&\quad:=p_1(u)+p_2(u)+p_3(u).
\EQNY
\nelem{lemP} leads to
\BQN\label{neg0}
p_i(u)=o\left(\frac{u}{m(u)\Delta_1(u)}\Psi(m(u))\right), \quad i=1,2,3
\EQN
implying that
\BQN\label{negl}
\Theta_0(u)=o\left(\frac{u}{m(u)\Delta_1(u)}\Psi(m(u))\right), \quad u\rw\IF.
\EQN
Since the proof \ehg{of} (\ref{neg0}) \he{is quite technical, we shall present it in} Appendix.
}
\subsubsection{Asymptotics of $\Theta(u)$}
\kd{
	We shall distinguish three scenarios:  $\varphi=0,$ $\varphi\in (0,\IF)$ and $\varphi=\IF$. The reason for this is that after rescaling the time of the correlation function in Lemma \ref{L2}, we get
\BQN\label{new1}
m^2(u)\left(1-r_u(\frac{\Delta_\gamma(u)s}{u},\frac{\Delta_1(u)t}{u},\frac{\Delta_\gamma(u)s_1}{u},\frac{\Delta_1(u)t_1}{u})\right)\sim \frac{\vf (\Delta_1(u)|t-t_1|)}{\vf(\Delta_1(u))}+\frac{\vf (\Delta_\gamma(u)|s-s_1|)}{\vf (\Delta_\gamma(u))}.
\EQN
If $\varphi=0$, then $\lim_{u\rw\IF}\Delta_\gamma(u)=0$ for $\gamma\in (0,1]$, implying that only the local behaviour of $\sigma^2$ at $0$ contributes to the limit in (\ref{new1}). If $\varphi\in (0,\IF)$, then $\lim_{u\rw\IF}\Delta_\gamma(u)\in (0,\IF)$, indicating that the whole function $\sigma^2$ determines the limit in (\ref{new1}). If $\varphi=\IF$, then  $\lim_{u\rw\IF}\Delta_\gamma(u)=\IF$, which means that the value of $\sigma^2(t)$  as $t\rw\IF$ is sufficient for the limit in (\ref{new1}).\\
\underline{\it \bf Case $\varphi=0$}:
We shall apply the {\it \he{uniform} double sum} technique  which is based
on appropriate division of the set $E(u)$ on "tiny" intervals for which one can uniformly derive exact asymptotics utilising our novel result in Lemma \ref{PICKANDS} in Appendix.
For this purpose we define}
$$F_{k,S}(u)=\Biggl[t_u+k\frac{\Delta_1(u)}{u}S, t_u+(k+1)\frac{\Delta_1(u)}{u}S
\Biggr],\ \ k\in \mathbb{Z}, S>0$$
$$L_{l,S}(u)=\Biggl[l\frac{\Delta_\gamma(u)}{u}S, (l+1)\frac{\Delta_\gamma(u)}{u}S
\Biggr],\ \ l\in \mathbb{N}\cup\{0\}, S>0$$
and set
\BQN\label{IK}
I_{k,l,S,S_1}(u)=L_{l,S_1}(u)\times F_{k,S}(u), \quad \eL{I_k(u):= I_{k,0,S,S_1} },
\EQN
 with $k\in \mathbb{N}, l\in \mathbb{N}\cup\{0\}, S, S_1>0 $.
Recall that, \kd{due to \eqref{Delta},
$\Delta_\gamma(u)=\overleftarrow{\sigma}
 \left(\frac{\sqrt{2}\sigma^2(u\eL{t_*})}{\gamma u(1+ct_*)}\right)$.
}
Further, let
 \BQN\label{NS}
 N_{S,u}=\left[\frac{\sigma(u)\ln u}{\Delta_1(u)S}\right]+1, \quad  N_{S_1,u}^{(1)}=\left[\frac{\overleftarrow{\sigma}\left(u^{-1}\sigma^2(u)\ln u\right)}{\Delta_\gamma(u)S_1}\right]+1
 \EQN
  and put
   $$\mathbb{V}_1=\{(k,k_1), -N_{S,u}\leq k<k_1\leq N_{S,u}, |k-k_1|>1 \},$$
   $$
   \mathbb{V}_2=\{(k,k_1), -N_{S,u}\leq k<k_1\leq N_{S,u}, k+1=k_1 \}.$$
\kd{We begin with the derivation of an upper estimate for $\Theta(u)$, as $u\to\infty$.}\\
 \underline{\it Upper bound of $\Theta(u)$}.
 Bonferroni inequality yields
\BQN\label{upper}
\Theta(u)&\leq&\sum_{k=-N_{S,u}}^{N_{S,u}}\pk{\sup_{(s,t)\in  \eL{I_k(u)} }Z_u(s,t)>m(u)}
\notag \\ &&+\sum_{k=-N_{S,u}}^{N_{S,u}}\sum_{l=1}^{N_{S_1,u}^{(1)}}\pk{\sup_{(s,t)\in I_{k,l,S,S_1}(u)}Z_u(s,t)>m(u)}\nonumber\\
&:=& \Theta_1(u)+\Theta_2(u).
\EQN
 In light of Lemma \ref{L1} for $u$ large enough
\BQN\label{maine1}
\Theta_1(u)\leq \sum_{k=-N_{S,u}}^{N_{S,u}}\pk{\sup_{(s,t)\in  \eL{I_k(u)} }\frac{\overline{Z}_u(s,t)}{1+(a_2-\epsilon)\frac{\sigma^2(us)}{\sigma^2(u)}}>m_{k,0}^{-\epsilon}(u)},
\EQN
with $\ve \in (0, \ehg{\min(a_1,a_2)})$ and  $$m_{k,0}^{\pm\epsilon}(u)=m(u)\left(1+(a_1-\epsilon)\left(k^*\frac{\Delta_1(u)}{u}S\right)^2\right),  \quad k^*=\min(|k|, |k+1|).$$
\eeh{\ehg{In order to derive an upper bound for}  $\Theta_1(u)$, we apply Lemma \ref{PICKANDS}
in Appendix, which gives uniform asymptotics for all terms in (\ref{maine1}). For this purpose, let  }
\rd{\BQN\label{TH11}
g_{u,k}=m_{k,0}^{-\epsilon}(u), \quad \xi_{u,k}=\frac{Z_{u,k}(s,t)}{1+f_{u,k}(s,t)}, \quad (s,t)\in \CE{E}=[0,S_1]\times [0,S],
 \EQN
with $\quad k\in K_u=\{k: -N_{S,u}\leq k \leq N_{S,u}\},$ where
$$Z_{u,k}(s,t)=\overline{Z}_u\left(\frac{\Delta_\gamma(u)}{u}s, t_{u,k}+\frac{\Delta_1(u)}{u}t\right), \quad  f_{u,k}(s)=(a_2-\epsilon)
\frac{\sigma^2(\Delta_\gamma(u)s)}{\sigma^2(u)}, \quad s\in [0,S_1]$$
and for $u>0$
$$ t_{u,k}=t_u+k\frac{\Delta_1(u)}{u}S.$$}
\kd{We check that
the conditions of
Lemma \ref{PICKANDS}
hold with the above introduced notation.
We start off with proving that
{\bf P1-P3} (see Appendix) hold with $$V(s,t)=B_{\alpha_0}(s)+B_{\alpha_0}^{\ehg{*}}(t), \quad
(s,t)\in [0,S_1]\times [0,S],$$
 where $B_{\alpha_0}$ and $B_{\alpha_0}^{\ehg{*}}$ are independent \eL{fBm's} with index $\alpha_0$.}
 It is straightforward that condition {\bf P1} holds. \rd{For {\bf P2}, by Lemma \ref{L2} and the fact that
	$$g_{u,k}\sim m(u), \quad u\rw\IF$$
	 uniformly with respect to $k\in K_u$, we have that for all $k\in K_u$ and $(s,t), (s_1,t_1)\in E$, as $u\rw\IF$
\BQNY\label{new0}
(g_{u,k})^2Var\left(Z_{u,k}(s,t)-Z_{u,k}(s_1,t_1)\right)&=&2(g_{u,k})^2\left(1-r_u\left(\frac{\Delta_\gamma(u)}{u}s,
\frac{\Delta_1(u)}{u}t,\frac{\Delta_\gamma(u)}{u}s_1, \frac{\Delta_1(u)}{u}t_1\right)\right)\nonumber\\
&\sim& (g_{u,k})^2\frac{\sigma^2(\Delta_1(u)|t-t_1|)+\gamma^2\sigma^2(\Delta_\gamma(u)|s-s_1|)}{\sigma^2(ut^*)}\nonumber\\
&\sim& 2\left(\frac{\sigma^2(\Delta_1(u)|t-t_1|)}{\sigma^2(\Delta_1(u))}+\frac{\sigma^2(\Delta_\gamma(u)|s-s_1|)}{\sigma^2(\Delta_\gamma(u))}\right),
\EQNY
implying that we can set
\BQN\label{theta}
\theta_{u,k}(s,t,s_1,t_1)=\frac{\sigma^2(\Delta_1(u)|t-t_1|)}{\sigma^2(\Delta_1(u))}+\frac{\sigma^2(\Delta_\gamma(u)|s-s_1|)}{\sigma^2(\Delta_\gamma(u))}, \quad (s,t),(s_1,t_1)\in \he{E}, k\in K_u.\EQN
 }
 Moreover, since
 $$\lim_{u\rw\IF}\Delta_\gamma(u)=0, \quad \gamma\in (0,1]$$
 by UCT
\BQNY
&&\lim_{u\rw\IF}\sup_{k\in K_u}\sup_{(s,t), (s_1,t_1)\in E}\left|\theta_{u,k}(s,t,s_1,t_1)-|s-s_1|^{2\alpha_0}-|t-t_1|^{2\alpha_0}\right|\\
&&= \lim_{u\rw\IF}\sup_{(s,t), (s_1,t_1)\in E}\left|\frac{\sigma^2(\Delta_\gamma(u)|s-s_1|)}{\sigma^2(\Delta_\gamma(u))}+\frac{\sigma^2(\Delta_1(u)|t-t_1|)}{\sigma^2(\Delta_1(u))}
-|s-s_1|^{2\alpha_\IF}-|t-t_1|^{2\alpha_0}\right|=0.
\EQNY
This means that {\bf P2} holds.
\COM{\KD{Since $\sigma^2$ is regularly varying at 0, then by
UCT \rd{and the fact that $g_{u,k}\sim m(u)$ uniformly with respect to $k\in K_u$,}
the second term on the right hand side of the above inequality satisfies}
\BQN\label{EQ1}
 &&\left|\frac{\sigma^2(\Delta_1(u)|t-t_1|)}{\sigma^2(\Delta_1(u))}\frac{(g_{u,k})^2\sigma^2(\Delta_1(u))}
{2\sigma^2(u\eL{t_*})}-|t-t_1|^{2\alpha_0}\right|\nonumber\\
&\quad&\quad \leq  \left|\frac{\sigma^2(\Delta_1(u)|t-t_1|)}{\sigma^2(\Delta_1(u))}-|t-t_1|^{2\alpha_0}\right|+\frac{\sigma^2(\Delta_1(u)|t-t_1|)}{\sigma^2(\Delta_1(u))}
\left|\frac{(g_{u,k})^2\sigma^2(\Delta_1(u))}
{2\sigma^2(u\eL{t_*})}-1\right|\rw 0
\EQN
uniformly with respect to $t,t_1\in [0,S]$ and $k\in K_u$. Similarly, we get that the first term also uniformly tends to 0 with respect to $s,s_1\in [0,S_1]$. Therefore we conclude that
$$\left|\theta_{u,k}(s,t,s_1,t_1)-|s-s_1|^{2\alpha_0}-|t-t_1|^{2\alpha_0}\right|\rw 0$$
 uniformly with respect to $(s,t), (s_1,t_1)\in E$ and $k\in K_u$, which implies that {\bf P2} holds.}

For {\bf P3}, by (\ref{re0}) we have that for $u$ sufficiently large
\BQNY
\theta_{u,k}(s,t,s_1,t_1))= \frac{\sigma^2(\Delta_1(u)|t-t_1|)}{\sigma^2(\Delta_1(u))}+\frac{\sigma^2(\Delta_\gamma(u)|s-s_1|)}{\sigma^2(\Delta_\gamma(u))}
\leq 2\left(S^{2\alpha_0-\lambda}+S_1^{2\alpha_0-\lambda}\right)\left(|s-s_1|^{\lambda}+|t-t_1|^{\lambda}\right)
\EQNY
for $(s,t), (s_1,t_1)\in \he{E}$ and all $k\in K_u$ with $0<\lambda<\min(2\alpha_0, 2\alpha_\IF)$.
 By UCT, we have for all $(s,t), (s_1,t_1)\in E$
\BQN
&&\sup_{|(s,t)-(s_1,t_1)|<\epsilon}|\theta_{u,k}(s,t,0,0)-\theta_{u,k}(s_1,t_1,0,0)|\nonumber\\
&\quad&\quad\leq
\sup_{|(s,t)-(s_1,t_1)|<\epsilon}\left|\frac{\sigma^2(\Delta_1(u)t)-\sigma^2(\Delta_1(u)t_1)}{\sigma^2(\Delta_1(u))}+\frac{\sigma^2(\Delta_\gamma(u)s)-
\sigma^2(\Delta_\gamma(u)s_1)}{\sigma^2(\Delta_\gamma(u))}-(t^{2\alpha_0}-t_1^{2\alpha_0}+s^{2\alpha_0}-s_1^{2\alpha_0})\right|\nonumber\\
&\quad&\quad\quad +\sup_{|(s,t)-(s_1,t_1)|<\epsilon}|t^{2\alpha_0}-t_1^{2\alpha_0}+s^{2\alpha_0}-s_1^{2\alpha_0}|\nonumber\\
&\quad&\quad \leq 2\epsilon+\sup_{|(s,t)-(s_1,t_1)|<\epsilon}|t^{2\alpha_0}-t_1^{2\alpha_0}+s^{2\alpha_0}-s_1^{2\alpha_0}|\leq \mathbb{C}\epsilon^{\alpha_0}, \ \ u\rw\IF,
\EQN
with $\mathbb{C}$ depending \ehc{only} on $\alpha_0$ (\eE{but not} on $k\in K_u$).
Moreover, using UCT, we have for $(s,t), (s_1,t_1)\in E$, $|(s,t)-(s_1,t_1)|<\epsilon$ and $k\in K_u$
\BQNY
\left|(g_{u,k})^2\left(1-r_u \Bigl(s_{u,l}+\frac{\Delta_\gamma(u)}{u}s, t_{u,k}+\frac{\Delta_1(u)}{u}t,s_{u,l}, t_{u,k}\Bigr)\right)-\theta_{u,k}(s,t,0,0)\right|
&\leq & \epsilon |\theta_{u,k}(s,t,0,0)|\\
& \leq& 2(S^{2\alpha_0}+S_1^{2\alpha_0})\epsilon
\EQNY
for all $u$ large. Consequently, as $u\to \IF$
\BQNY
&&(g_{u,k})^2\mathbb{E}\{\left[Z_{u,k}(s,t)-Z_{u,k}(s_1,t_1)\right]Z_{u,k}(0,0)\}\\
&\quad&\quad \leq \left|(g_{u,k})^2\left(1-r_u(s_{u,l}+\frac{\Delta_\gamma(u)}{u}s, t_{u,k}+\frac{\Delta_1(u)}{u}t,s_{u,l}, t_{u,k})\right)-\theta_{u,k}(s,t,0,0)\right| \\
&\quad&\quad +\left|(g_{u,k})^2\left(1-r_u(s_{u,l}+\frac{\Delta_\gamma(u)}{u}s_1, t_{u,k}+\frac{\Delta_1(u)}{u}t_1,s_{u,l}, t_{u,k})\right)-\theta_{u,k}(s,t,0,0)\right|\\
&\quad&\quad +|\theta_{u,k}(s,t,0,0)-\theta_{u,k}(s,t,0,0)|\\
&\quad&\quad \leq\mathbb{C}\epsilon^{\alpha_0}+4(S^{2\alpha_0}+S_1^{2\alpha_0})\epsilon
\EQNY
 uniformly for $(s,t), (s_1,t_1)\in E$, $|(s,t)-(s_1,t_1)|<\epsilon$ and $k\in K_u$.
Letting $\epsilon\rw 0$, we confirm that ${\bf P3}$ holds. Hence
we can conclude that {\bf P1-P3} hold with $V(s,t)=B_{\alpha_0}(s)+B_{\alpha_0}^{\ehg{*}}(t), (s,t)\in E$, where $B_{\alpha_0}$ and $B_{\alpha_0}^{\ehg{*}}$ are independent \eL{fBm's} with index $\alpha_0$. Therefore, by the fact that for all $\ve>0$ sufficiently small \ehb{(hereafter $\RW$ means uniform convergence)}
$$g_{u,k}^2f_{u,k}(s)\RW \gamma_\epsilon s^{2\alpha_0}, s\in[0,S_1], \quad \text{with} \quad \gamma_\epsilon=\frac{a_2-\epsilon}{a_2} \barg,$$
 and Lemma \ref{PICKANDS}
\KD{we have}
 \BQN\label{TH12}
\frac{\pk{\sup_{(s,t)\in  E }\xi_{u,k}(s,t)>g_{u,k}}}{\Psi(g_{u,k})}\rw \mathcal{R}_V^{ \gamma_\epsilon s^{2\alpha_0}}(E)=\mathcal{H}_{B_{\alpha_0}}[0,S]\mathcal{P}_{B_{\alpha_0}}^{\gamma_\epsilon}[0,S_1], \ \ u\rw\IF
 \EQN
 uniformly with respect to $-N_{S,u}\leq k \leq N_{S,u}$.
From (\ref{maine1}) and (\ref{TH12}) it follows that
\BQN\label{pi1}
\Theta_1(u)&\leq& \sum_{k=-N_{S,u}}^{N_{S,u}}\mathcal{H}_{B_{\alpha_0}}[0,S]\mathcal{P}_{B_{\alpha_0}}^{\gamma_\epsilon}[0,S_1]\Psi(m_{k,0}^{-\epsilon}(u))(1+o(1))\nonumber\\
&\leq&\sum_{k=-N_{S,u}}^{N_{S,u}}\mathcal{H}_{B_{\alpha_0}}[0,S]\mathcal{P}_{B_{\alpha_0}}^{\gamma_\epsilon}[0,S_1]\Psi(m(u))e^{-(a_1-\epsilon)
\left(k^*m(u)\frac{\Delta_1(u)}{u}S\right)^2}(1+o(1))\nonumber\\
&\leq&\frac{\mathcal{H}_{B_{\alpha_0}}[0,S]}{S}\mathcal{P}_{B_{\alpha_0}}^{\gamma_\epsilon}[0,S_1](a_1-\epsilon)^{-1/2}
\Psi(m(u))\frac{u}{m(u)\Delta_1(u)}\int_{-\IF}^{\IF}e^{-x^2}dx(1+o(1))\nonumber\\
&\sim&(a_1-\epsilon)^{-1/2}\sqrt{\pi}\mathcal{H}_{B_{\alpha_0}}\mathcal{P}_{B_{\alpha_0}}^{\gamma_\epsilon}\Psi(m(u))\frac{u}{m(u)\Delta_1(u)} ,\ \ \EQN
\eeh{as $u\rw\IF, S, S_1\rw\IF$ (in this order).}\\
Next, we deal with  $\Theta_2(u)$. By UCT, for any $\epsilon>0$ sufficiently small
\BQNY
\he{(a_2-\epsilon)} \sup_{s\in E_1(u) }\frac{\sigma^2(us)}{\sigma^2(u)}\rw 0, \quad u\rw\IF.
\EQNY
Moreover, by (\ref{re00}) for $u$ large enough
 \BQNY
\inf_{s\in L_{l,S}(u)}\left(m_{k,0}^{-\epsilon}(u)\right)^2\frac{\sigma^2(us)}{\sigma^2(u)}&\geq &  \frac{1}{2}\inf_{s\in [lS_1, (l+1)S_1]}\frac{\sigma^2(\Delta_\gamma(u)s)}{\sigma^2(\Delta_\gamma(u))}\frac{\sigma^2(\Delta_\gamma(u))}{\sigma^2(u)}m^2(u)\nonumber\\
&\geq &  Q\inf_{s\in [lS_1, (l+1)S_1]}\frac{\sigma^2(\Delta_\gamma(u)s)}{\sigma^2(\Delta_\gamma(u))}\nonumber\\
&\geq&  Q(lS_1)^{{\lambda}},\ \ 1 \leq l\leq N_{S_1,u}^{(1)}, 0<\lambda<\min (2\lambda_0, 2\alpha_\IF).
\EQNY
Consequently, by Lemma \ref{L1} and Lemma \ref{PICKANDS} \EH{(note that we can similarly show  the validity of  {\bf P1-P3} for \ehg{$\overline{Z}_u(s,t)$})}  we have  for any $\epsilon>0$
\BQN\label{pi2}
\Theta_2(u)&\leq& \sum_{k=-N_{S,u}}^{N_{S,u}}\sum_{l=1}^{N_{S_1,u}^{(1)}}\pk{\sup_{(s,t)\in I_{k,l,S,S_1}(u)}\overline{Z}_u(s,t)>m_{k,0}^{-\epsilon}(u)\left(1+(a_2-\epsilon)\inf_{s\in L_{l,S}(u) }\frac{\sigma^2(us)}{\sigma^2(u)}\right)}\nonumber\\
&\leq&\sum_{k=-N_{S,u}}^{N_{S,u}}\sum_{l=1}^{N_{S_1,u}^{(1)}}\mathcal{H}_{B_{\alpha_0}}[0,S]\mathcal{H}_{B_{\alpha_0}}[0,S_1]
\Psi(m_{k,0}^{-\epsilon}(u))e^{-Q_1(lS_1)^{{\lambda}}}(1+o(1))\nonumber\\
&\leq& \sum_{k=-N_{S,u}}^{N_{S,u}}\mathcal{H}_{B_{\alpha_0}}[0,S]\mathcal{H}_{B_{\alpha_0}}[0,S_1]
\Psi(m_{k,0}^{-\epsilon}(u))e^{-Q_2 S_1^{ \CE{\lambda} }}(1+o(1))\nonumber\\
&=&o\left(\Psi(m(u))\frac{u}{m(u)\Delta_1(u)}
\right), \ \
\EQN
\eeh{as $u\rw\IF,  S,S_1\rw\IF$.}
Combining (\ref{pi1}) and (\ref{pi2}), and letting $\epsilon\rw 0$, we derive the upper bound of $\Theta(u)$. \\
\underline{\it Lower bound of $\Theta(u)$.} By Bonferroni inequality we obtain
\BQN\label{lower}
\Theta(u)\geq \sum_{k=-N_{S,u}+1}^{N_{S,u}-1}\pk{\sup_{(s,t)\in  \eL{I_k(u)} }Z_u(s,t)>m(u)}-\Sigma_1(u)-\Sigma_2(u):=J(u)-\Sigma_1(u)-\Sigma_2(u),
\EQN
with
\BQN\label{SIG}
\Sigma_i(u)=\sum_{(k,k_1)\in \mathbb{V}_i}\pk{\sup_{(s,t)\in  \eL{I_k(u)} }Z_u(s,t)>m(u), \sup_{(s,t)\in  \eL{I_{k_1}(u)} }Z_u(s,t)>m(u)}, \ \ i=1,2.
\EQN
With similar arguments as in the \EH{derivation} of (\ref{pi1}) we have
\BQN\label{LOW1}
J(u)\geq (a_1+\epsilon)^{-1/2}\sqrt{\pi}\mathcal{H}_{B_{\alpha_0}}\mathcal{P}_{B_{\alpha_0}}^{\gamma_{-\epsilon}}\Psi(m(u))\frac{u}{m(u)\Delta_1(u)}(1+o(1)),\ \ u\rw\IF, S, S_1\rw\IF.
\EQN
In light of Lemma \ref{L2} and (\ref{re00})
we have for $(s,t,s_1,t_1)\in  \eL{I_k(u)} \times  \eL{I_{k_1}(u)} $ with $(k,k_1)\in \mathbb{V}_1$
\def\ldy{{\kappa}}
\BQNY
2\leq \eL{Var}(\overline{Z}_u(s,t)+\overline{Z}_u(s_1,t_1))&=&4-2(1-r_u(s,t,s_1,t_1))\\&\leq& 4-\frac{\gamma^2\sigma^2(u|s-s_1|)+\sigma^2(u|t-t_1|)}{2\sigma^2(u\eL{t_*})}\\
&\leq& 4-\frac{1}{2m^2(u)}
\frac{\sigma^2(\Delta_1(u)|u(t-t_1)/\Delta_1(u)|)}{\sigma^2(\Delta_1(u))}\\
 &\leq& 4-Q_3\frac{|k_1-k|^{ \CE{\lambda} }S^{ \CE{\lambda} }}{m^2(u)},
\EQNY
where $0< \CE{\lambda} < \min(2\alpha_0, 2\alpha_\IF)$, implying thus
\BQNY
\Sigma_1(u)&\leq& \sum_{(k,k_1)\in \mathbb{V}_1}\pk{\sup_{(s,t)\in  \eL{I_k(u)} }\overline{Z}_u(s,t)>m_{k,0}^{-\epsilon}(u), \sup_{(s,t)\in  \eL{I_{k_1}(u)} }\overline{Z}_u(s,t)>m_{k_1,0}^{-\epsilon}(u)}\\
&\leq&\sum_{(k,k_1)\in \mathbb{V}_1}\pk{\sup_{(s,t,s_1,t_1)\in  \eL{I_k(u)} \times  \eL{I_{k_1}(u)} }
\Bigl(\overline{Z}_u(s,t)+\overline{Z}_u(s_1,t_1)\Bigr)>2\widetilde{m}_{k,k_1,0}^{-\epsilon}(u)}\\
&\leq&\sum_{(k,k_1)\in \mathbb{V}_1}\pk{\sup_{(s,t,s_1,t_1)\in  \eL{I_k(u)} \times  \eL{I_{k_1}(u)} }\Bigl(\overline{\overline{Z}_u(s,t)+\overline{Z}_u(s_1,t_1)}\Bigr)>\frac{2\widetilde{m}_{k,k_1,0}^{-\epsilon}(u)}{\sqrt{4-Q_3
\frac{|k_1-k|^{ \CE{\lambda} }S^{ \CE{\lambda} }}{m^2(u)}}}},
\EQNY
with $\widetilde{m}_{k,k_1,0}^{-\epsilon}(u)=\min(m_{k,0}^{-\epsilon}(u),m_{k_1,0}^{-\epsilon}(u))$.\\
Let next
$$r_u(t,s,t_1,s_1,t',s',t_1^{'},s_1^{'})=\eL{Cor}(\overline{Z_u}(s,t)+\overline{Z_u}(s_1,t_1),
\overline{Z_u}(s',t')+\overline{Z_u}(s_1^{'},t_1^{'})).$$  \eL{By} Lemma \ref{L2} and (\ref{re00}), for $u$ sufficiently large
\BQNY
1-\eL{r_u(s,t,s_1,t_1,s',t',s_1',t_1^{'})}&\leq& \frac{\eL{Var}( \overline{Z_u}(s,t)+\overline{Z_u}(s_1,t_1)-\overline{Z_u}(s',t')-\overline{Z_u}(s_1^{'},t_1^{'}))}
{2\sqrt{\eL{Var}(\overline{Z_u}(s,t)+\overline{Z_u}(s_1,t_1))}\sqrt{\eL{Var}(\overline{Z_u}(s',t')+\overline{Z_u}(s_1^{'},t_1^{'})}}\\
&\leq& 1-r_u(s,t,s',t')+1-r_u(s_1,t_1,s_1^{'},t_1^{'})\\
&\leq& \frac{2}{m^2(u)}
\frac{\sigma^2(\Delta_\gamma(u)|u(s-s_1)/\Delta_\gamma(u)|)+\sigma^2(\Delta_\gamma(u)|u(s'-s_1')/\Delta_\gamma(u)|)}{\sigma^2(\Delta_\gamma(u))}\\
&&+\frac{2}{m^2(u)}
\frac{\sigma^2(\Delta_1(u)|u(t-t_1)/\Delta_1(u)|)+\sigma^2(\Delta_1(u)|u(t'-t_1')/\Delta_1(u)|)}{\sigma^2(\Delta_1(u))}\\
&\leq&\frac{ Q_4(S^*)^2}{m^2(u)} \Biggl[\left(\frac{u}{\Delta_\gamma(u)}\right)^{\ldy}
\left(|s-s'|^{\ldy}+|s_1-s_1^{'}|^{\ldy}\right)+
\left(\frac{u}{\Delta_1(u)}\right)^{\ldy}\left(|t-t'|^{\ldy}+|t_1-t_1^{'}|^{\ldy}\right)\Biggr]
\EQNY
holds for $(t,s,t_1,s_1)$, $(t',s',t_1^{'},s_1^{'})$ $\in  \eL{I_k(u)} \times  \eL{I_{k_1}(u)} $ with $0<\ldy<\min(2{\alpha_\IF}, 2\alpha_0)$ and $S^*=\max(S,S_1)\geq 1$.
Define the \blu{following} homogeneous Gaussian field
$$X_u^*(s,t,s_1,t_1)=2^{-1}(X^{1}_u(s)+X^{2}_u(t)+X^{3}_u(s_1)+X^{4}_u(t_1)), \quad (s,t,s_1,t_1)\in \mathbb{R}^4,$$
 with $X_u^{i}(s), 1\leq i\leq 4,$ being independent \eL{with} the correlation functions
\BQNY
r_u^{(i)}(s,s')=e^{-8Q_4(S^*)^2\left(\frac{u}{\Delta_1(u)}\right)^{\ldy}\frac{
|s-s'|^{\ldy}}{m^2(u)}}, \ \ i=1,3,
\EQNY
\BQNY
r_u^{(i)}(s,s')=e^{-8Q_4(S^*)^2\left(\frac{u}{\Delta_\gamma(u)}\right)^{\ldy}\frac{
|s-s'|^{\ldy}}{m^2(u)}},\ \ i=2,4.
\EQNY
We denote the correlation function of $X_u^*$ by $r_u^{*}$.
Clearly, for $(t,s,t_1,s_1)$, $(t',s',t_1^{'},s_1^{'})$ $\in  \eL{I_k(u)} \times  \eL{I_{k_1}(u)} $ and $u$ large enough
\BQNY
1-r_u(s,t,s_1,t_1,s',t',s_1^{'},t_1^{'})
\leq 1-r_u^{*}(s,t,s_1,t_1,s',t',s_1^{'},t_1^{'}).
\EQNY
In light of Slepian's inequality (see e.g., \eG{Theorem 2.2.1 in} \cite{AdlerTaylor}\EH{; \ehc{note in passing that there is a remarkable extension of this inequality for stable processes, see \cite{GennaSlepian}}}) and Lemma \ref{PICKANDS} we have
\BQN\label{LOW2}
\Sigma_1(u)&\leq& \sum_{(k,k_1)\in \mathbb{V}_1}\pk{\sup_{(s,t,s_1,t_1)\in  \eL{I_k(u)} \times  \eL{I_{k_1}(u)} }X_u^*(s,t,s_1,t_1)>\frac{2\widetilde{m}_{k,k_1,0}^{-\epsilon}(u)}{\sqrt{4-Q_3
\frac{|k_1-k|^{ \CE{\lambda} }S^{ \CE{\lambda} }}{m^2(u)}}}}\nonumber\\
&\leq&\sum_{(k,k_1)\in \mathbb{V}_1} \left(\mathcal{H}_{B_{\ldy/2}}[0,S_2]\right)^2\left(\mathcal{H}_{B_{\ldy/2}}[0,S_3]\right)^2
\Psi\left(\frac{2\widetilde{m}_{k,k_1,0}^{-\epsilon}(u)}{\sqrt{4-Q_3
\frac{|k_1-k|^{ \CE{\lambda} }S^{ \CE{\lambda} }}{m^2(u)}}}\right)(1+o(1))\nonumber\\
&\leq& \sum_{k=-N_{S,u}}^{N_{S,u}}\left(\frac{\mathcal{H}_{B_{\ldy/2}}[0,S_2]}{S_2}\right)^2
\left(\frac{\mathcal{H}_{B_{\ldy/2}}[0,S_3]}{S_3}\right)^2\Psi(m_{k,0}^{-\epsilon}(u))S_2^{-2}S_3^{-2}\sum_{j\geq 1}e^{-Q_5(jS)^{ \CE{\lambda} }}(1+o(1))\nonumber\\
&\leq & Q_6\Psi(m(u))\frac{u}{m(u)\Delta_1(u)}S_1^{-2}e^{-Q_7S^{ \CE{\lambda} }}(1+o(1))\notag \\
&=&o\left(\Psi(m(u))\frac{u}{m(u)\Delta_1(u)}\right),\ \
\EQN
with $S_2=(2Q_4(S^*)^2)^{1/\ldy}S$ and $S_3=(2Q_4(S^*)^2)^{1/\ldy}S_1$, \eeh{as $u\rw\IF, S\rw\IF$
	\he(in this order)}. Further, we obtain
\BQN\label{LOW3}
\Sigma_2(u)&=&\sum_{k=-N_{S,u}}^{N_{S,u}}\pk{\sup_{(s,t)\in  \eL{I_k(u)} }Z_u(s,t)>m(u), \sup_{(s,t)\in \eL{I_{k+1}(u)}}Z_u(s,t)>m(u)}\nonumber\\
&\leq&\sum_{k=-N_{S,u}}^{N_{S,u}}\left[\pk{\sup_{(s,t)\in  \eL{I_k(u)} }Z_u(s,t)>m(u)}+\pk{ \sup_{(s,t)\in \eL{I_{k+1}(u)}}Z_u(s,t)>m(u)}\right.\nonumber\\
&&\left.-\pk{\sup_{(s,t)\in \eL{I_{k}(u)}\cup \eL{I_{k+1}(u)}}Z_u(s,t)>m(u)}\right]\nonumber\\
&\leq& \sum_{k=-N_{S,u}}^{N_{S,u}}\left(2\mathcal{H}_{B_{\alpha_0}}[0,S]-\mathcal{H}_{B_{\alpha_0}}[0,2S]\right)
\mathcal{P}_{B_{\alpha_0}}^{\gamma_\epsilon}[0,S_1]\Psi(\widetilde{m}_{k,k+1,0}^{-\epsilon}(u))(1+o(1))\nonumber\\
&=&o\left(\Psi(m(u))\frac{u}{m(u)\Delta_1(u)}\right)
\EQN
\eeh{ as $ u\rw\IF, S_1\rw\IF, S\rw\IF$.}
Combining (\ref{LOW1})-(\ref{LOW3})
and letting $\epsilon\rw0$, the lower bound of $\Theta(u)$ is derived.
Since the upper and lower bound \KD{coincide}, then we have
\BQNY
\Theta(u)\sim \sqrt{\frac \pi {a_1}}\mathcal{H}_{B_{\alpha_0}}\mathcal{P}_{B_{\alpha_0}}^{ \barg }\Psi(m(u))\frac{u}{m(u)\Delta_1(u)}
\EQNY
\eE{and therefore the claim follows  by (\ref{P1}) and (\ref{negl})}- \\
 \underline{\it \bf Case $\varphi\in (0,\IF)$}:  \ehd{The main difference to the above proof is that
$\Delta_\gamma(u)=1$ and $\gamma\in (0,1]$, which influences  (\ref{theta}) and (\ref{TH12}) and hence the
resulting Pickands or Piterbarg constants that show up in the result}. Therefore, in order to avoid repetitions,
we present only the counterpart of the derivations of (\ref{theta}) and (\ref{TH12}).
Next, we check {\bf P2-P3} (\CE{conditions {\bf P1} is  easy to verify}) by using the same notation as in (\ref{TH11}) and (\ref{theta}).
In order to prove {\bf P2}, in light of Lemma \ref{L2} and the fact that $g_{u,k}\sim m(u)$ as $u\rw\IF$ uniformly with respect to $k\in K_u$, we have that for all $k\in K_u$ and $(s,t), (s_1,t_1)\in E$, as $u\rw\IF$
\BQNY
(g_{u,k})^2Var\left(Z_{u,k}(s,t)-Z_{u,k}(s_1,t_1)\right)&=&2(g_{u,k})^2\left(1-r_u\left(\frac{\Delta_\gamma(u)}{u}s,
\frac{\Delta_1(u)}{u}t,\frac{\Delta_\gamma(u)}{u}s_1, \frac{\Delta_1(u)}{u}t_1\right)\right)\\
&\sim& (m(u))^2\frac{\sigma^2(\Delta_1(u)|t-t_1|)+\gamma^2\sigma^2(\Delta_\gamma(u)|s-s_1|)}{\sigma^2(ut^*)}\\
&\sim& 2\left(\frac{2 c^2\gamma^2}{\varphi^2}\sigma^2(|s-s_1)+\frac{2 c^2}{\varphi^2}\sigma^2(|t-t_1|)\right).
\EQNY
Hence, we can set that
 \BQNY
\theta_{u,k}(s,t,s_1,t_1)=\frac{2 c^2\gamma^2}{\varphi^2}\sigma^2(|s-s_1)+\frac{2 c^2}{\varphi^2}\sigma^2(|t-t_1|), \quad (s,t), (s_1,t_1)\in E, k\in K_u,
 \EQNY
  which ensures that {\bf P2} holds. Next, for {\bf P3}, in light of (\ref{re00})  we derive that for $u$ sufficiently large and $\lambda\in (0,\min (2\alpha_0, 2\alpha_\IF))$,
 \BQNY
\theta_{u,k}(s,t,s_1,t_1)
 = \frac{2c^2}{\varphi^2}|\sigma^2(|s-s_1|)+\sigma^2(|t-t_1|)|
 \leq Q\left(|s-s_1|^{ \CE{\lambda} }+t-t_1|^{ \CE{\lambda} }\right),
 \EQNY
with $k\in K_u, (s,t),(s_1,t_1)\in E$.
 In addition, for $(s,t), (s_1,t_1)\in E$, $|(s,t)-(s_1,t_1)|<\epsilon$, $k\in K_u$ and $u$ sufficiently large we have
 \BQNY
 &&(g_{u,k})^2\mathbb{E}\{\left[Z_{u,k,l}(s,t)-Z_{u,k,l}(s_1,t_1)\right]Z_{u,k,l}(0,0)\}\\
 &\quad&\quad\leq \left|(g_{u,k})^2\left(1-r_u(s_{u,l}+\frac{1}{u}s, t_{u,k}+\frac{1}{u}t,s_{u,l}, t_{u,k})\right)-\theta_{u,k}(s,t,0,0)\right| \\
 &\quad&\quad \ \ \ +\left|(g_{u,k})^2\left(1-r_u(s_{u,l}+\frac{1}{u}s_1, t_{u,k}+\frac{1}{u}t_1,s_{u,l}, t_{u,k})\right)-\theta_{u,k}(s_1,t_1,0,0)\right|\\
 &\quad&\quad \ \ \ +|\theta_{u,k}(s,t,0,0)-\theta_{u,k}(s_1,t_1,0,0)|\\
 &\quad&\quad\leq\epsilon |\theta_{u,k}(s,t,0,0)+\theta_{u,k}(s_1,t_1,0,0)|+|\theta_{u,k}(s,t,0,0)-\theta_{u,k}(s_1,t_1,0,0)|\\
 &\quad&\quad\leq\mathbb{C}_2\left(\epsilon+\left|\sigma^2(t)-\sigma^2(t_1)\right|+\left|\sigma^2(s)-\sigma^2(s_1)\right|\right)\rw 0, \quad \epsilon\rw 0.
 \EQNY
 Thus {\bf P3} is satisfied. Next let
$$V(s,t)=\frac{1+c\eL{t_*}}{\sqrt{2}\varphi \eL{t_*}^{2\alpha_\IF}} \eHH{\bigl[} \gamma X(s)+X^{\ehg{*}}(t)\bigr]=\frac{\sqrt{2}c}{\varphi } \eHH{\bigl[} \gamma X(s)+X^{\ehg{*}}(t)\bigr], \quad  (s,t)\in E,$$ with $X^{\ehg{*}}$ an independent copy of $X$. Hence by Lemma \ref{PICKANDS} and the fact that (recall that  $\gamma_\epsilon=\frac{a_2-\epsilon}{a_2} \barg$)
$$(g_{u,k})^2f_{u,k}(s,t)\RW \frac{\gamma_\epsilon\gamma^2(1+c\eL{t_*})^2}{2\eL{t_*}^{4\alpha_\IF}\varphi^2}\sigma^2(s)=
\frac{2\gamma_\epsilon c^2\gamma^2}{\varphi^2}\sigma^2(s), \quad (s,t)\in E, \ u\rw\IF,$$
we have
 \BQNY
 \frac{\pk{\sup_{(s,t)\in  E} \xi_{u,k}(s,t)>g_{u,k}}}{\Psi(g_{u,k})}
 \rw\mathcal{R}_V^{ \frac{2\gamma_\epsilon c^2\gamma^2}{\varphi^2}\sigma^2(s)}(E)
 =\mathcal{H}_{\frac{\sqrt{2}c}{\varphi }X}[0,S]\mathcal{P}_{\frac{\sqrt{2}c\gamma}{\varphi }X}^{\gamma_\epsilon}[0,S_1],\ \ u\rw\IF
 \EQNY
 uniformly  with respect to $k\in K_u$. Repeating the \ehd{derivations} of (\ref{pi1})-(\ref{LOW3}),
 we conclude that the claim follows with the generalised Pickands and Piterbarg constants above instead of those for case
 $\varphi=0$.  Note that the existence of $\mathcal{H}_{X^*}$
 has been proved, see e.g. \cite{Pit96}, \cite{DE2002} and \cite{DI2005}; the
 \ehd{proof of the} finiteness of the generalised Piterbarg constants $\lim_{S_1\rw\IF}\mathcal{P}_{\frac{\sqrt{2}c\gamma}{\varphi }X}^{\gamma_\epsilon}[0,S_1]$ is postponed to Lemma \ref{EXISTENCE} in the Appendix.\\
\rd{\underline{\it\bf Case $\varphi=\IF$}: Since $\Delta_\gamma(u)$ is the same as in the case $\varphi=0$, the proof is very similar to that case.
The main difference is that the limiting Gaussian process $V$  that appears in (\ref{TH12}) is here different, namely
{\bf P1-P3} hold with
$$V(s,t)=B_{\KD{\alpha_\IF}}(s)+B_{\KD{\alpha_\IF}}^{*}(t), \quad (s,t)\in [0,S_1]\times [0,S],$$
where $B_{\alpha_\IF}$ and $B_{\alpha_\IF}^{*}$ are independent \eL{fBm's} with index $\alpha_\IF$.
We omit details.}
\QED

\subsection{\prooftheo{THA}} \eeh{We begin with transformation of  the distribution of the conditional passage time
to the ratio of two tail probabilities of supremum of $\gamma$-reflected Gaussian process  over appropriately chosen intervals.}

\kd{Using the same notation as introduced} in the proof of Theorem \ref{TH1},
\eL{first we focus} on   $\tau_1^*(u)$. Let $D_{x,u}=\{(s,t): 0\leq s\leq t\leq xu^{-1}A(u)+t_u\}$. For all $u$ large we have
\BQN\label{RT1}
\pk{\frac{\tau_1^*(u)-ut_u}{A(u)}\leq x}&=&\frac{\pk{\tau_1(u)\leq x A(u)+ut_u}}{\pk{\tau_1(u)<\IF}}=\frac{\pk{\sup_{t\in[0,x A(u)+ut_u]}W_{\gamma}(t)>u}}{\psi_{\gamma,\IF}(u)}\nonumber\\
&=&\frac{\pk{\sup_{(s,t)\in D_{x,u}}Z_u(s,t)>m(u)}}{\psi_{\gamma,\IF}(u)},
\EQN
with $Z_u(s,t)$ defined in (\ref{ZU}) and $m(u)$ defined in (\ref{MU}).
\kd{By Theorem \ref{TH1}, it suffices to find
\kk{the asymptotics of \\
$\pk{\sup_{(s,t)\in D_{x,u}}Z_u(s,t)>m(u)}$,}
for which we write}
\BQN\label{RT2}
\pi^x(u)\leq \pk{\sup_{(s,t)\in D_{x,u}}Z_u(s,t)>m(u)}\leq \pi^x(u)+\pk{\sup_{(s,t)\in \eE{D \setminus } E(u)}Z_u(s,t)>m(u)},
\EQN
where $$\pi^x(u)=\pk{\sup_{(s,t)\in E_1(u)\times E_2^x(u)}Z_u(s,t)>m(u)}, \quad  E_2^x(u)=\left(t_u-\frac{\sigma(u)\ln u}{u}, t_u+xu^{-1}A(u)\right)$$
with $D$ defined in (\ref{ZU}) and  $E_1(u), E(u)$ defined in (\ref{EU}).
Moreover,
\BQN\label{RT3}
J^x(u)-\Sigma_1(u)-\Sigma_2(u)\leq \pi^x(u)\leq \pi_1^x(u)+\Theta_2(u),
\EQN
where
$$\pi_1^x(u)=\sum_{k=-N_{S,u}}^{N_{S,u}^x}\pk{\sup_{(s,t)\in  \eL{I_k(u)} }Z_u(s,t)>m(u)}, \ \ J^x(u)=\sum_{k=-N_{S,u}+1}^{N_{S,u}^x-1}\pk{\sup_{(s,t)\in  \eL{I_k(u)} }Z_u(s,t)>m(u)} $$
with
\kk{$N_{S,u}^x=\left[\frac{x A(u)}{\Delta_1(u)S}\right]+1$,}
$I_{k}(u)$ defined in (\ref{IK}), $N_{S,u}$ in (\ref{NS}), $\Theta_2(u)$ in (\ref{upper}) and  $\Sigma_i(u), i=1,2$ in (\ref{SIG}).\\
\underline{\it \bf Case $\varphi=0$}:
Similarly as in (\ref{pi1}), with $\epsilon\in (0,a_1)$ and $k^*=\min(|k|, |k+1|)$, we have that
\BQN\label{RT4}
\pi_1^x(u)&\leq&\sum^{N_{S,u}^x}_{k=-N_{S,u}}\mathcal{H}_{B_{\alpha_0}}[0,S]\mathcal{P}_{B_{\alpha_0}}^{\gamma_\epsilon}[0,S_1]\Psi(m(u))e^{-(a_1-\epsilon)
\left(k^*m(u)\frac{\Delta_1(u)}{u}S\right)^2}(1+o(1))\nonumber\\
&=&\frac{\mathcal{H}_{B_{\alpha_0}}[0,S]}{S}\mathcal{P}_{B_{\alpha_0}}^{\gamma_\epsilon}[0,S_1](a_1-\epsilon)^{-1/2}
\Psi(m(u))\frac{u}{m(u)\Delta_1(u)}\int_{-\IF}^{\sqrt{\frac{a_1-\epsilon}{2a_1}}x}e^{-y^2}dy(1+o(1))\nonumber\\
&\sim& \sqrt{\pi/a_1}\Phi(x)\mathcal{H}_{B_{\alpha_0}}\mathcal{P}_{B_{\alpha_0}}^{ \barg }\Psi(m(u))\frac{u}{m(u)\Delta_1(u)}, \ \
\EQN
\eeh{as $u\rw\IF, S,S_1\rw\IF, \epsilon\rw 0$}, and
\BQN\label{RT5}
J^x(u)\geq \sqrt{\pi/a_1}\Phi(x)\mathcal{H}_{B_{\alpha_0}}\mathcal{P}_{B_{\alpha_0}}^{ \barg }\Psi(m(u))\frac{u}{m(u)\Delta_1(u)}(1+o(1)).
\EQN
By (\ref{negl})
$$\pk{\sup_{(s,t)\in \eE{D \setminus } E(u)}Z_u(s,t)>m(u)}=o(\pi_1^x(u))=o(J^x(u)).$$
Furthermore, it follows from (\ref{pi2}), (\ref{LOW2}) and (\ref{LOW3}) that
$\Theta_2(u)$, $\Sigma_1(u)$ and $\Sigma_2(u)$ are all
negligible in comparison with $\pi_1^x(u)$ and $J^x(u)$ for $x\in (-\IF,\IF]$.
Therefore, as $u\to \IF$,
\BQN\label{RT6}
\pk{\sup_{(s,t)\in D_{x,u}}Z_u(s,t)>m(u)}\sim \sqrt{\pi/a_1} \Phi(x)\mathcal{H}_{B_{\alpha_0}}\mathcal{P}_{B_{\alpha_0}}^{\gamma_\epsilon}\Psi(m(u))\frac{u}
{m(u)\Delta_1(u)}\sim\Phi(x)\psi_{\gamma,\IF}(u),
\EQN
which together  with (\ref{RT1}) implies
\BQNY
\limit{u} \pk{\frac{\tau_1^*(u)-ut_u}{A(u)}\leq x}=\Phi(x),\ \ x\in(-\IF,\IF].
\EQNY
\kd{Next,} we investigate the last passage time. Similarly as above, for $x\in (-\IF,\IF]$ we have
\BQN\label{RT7}
\pk{\frac{\tau_2^*(u)-ut_u}{A(u)}\leq x}&=&1-\pk{\frac{\tau_2(u)-ut_u}{A(u)}> x\Bigl\lvert \tau_1(u)<\IF }\\
&=&1-\frac{\pk{\sup_{t\in [x A(u)+ut_u,\IF)}W_\gamma(t)>u}}{\pk{\tau_1(u)<\IF}}\nonumber\\
&=&1-\frac{\pk{\sup_{t\in [x u^{-1}A(u)+t_u,\IF)}Z_u(s,t)>m(u)}}{\pk{\tau_1(u)<\IF}}\nonumber\\
&\rw& 1-\Psi(x)=\Phi(x) \nonumber
\EQN
 as $u\rw\IF$. Hence application of Lemma 2.1 in \cite{HJ13} (recall that $\tau_1(u)\leq \tau_2(u)$) yields that for any $x, y\in \mathbb{R}$
 $$\pk{\frac{\tau_1^*-ut_u}{A(u)}\leq x, \frac{\tau_2^*-ut_u}{A(u)}\leq y}\rw \pk{\mathcal{N}\leq \min(x,y)}, \quad u\rw\IF.$$\\
 {\underline {\it \bf Case $\eL{\varphi} \in (0,\IF)$}:
Note that  (\ref{RT4}) and (\ref{RT5}) are also valid by replacing
 $\mathcal{H}_{B_{\alpha_0}}[0,S]$ with $\mathcal{H}_{\frac{\sqrt{2}c}{\varphi}X}[0,S]$ and
 $\mathcal{P}_{B_{\alpha_0}}^{\gamma_\epsilon}[0,S_1]$ with
 $\mathcal{P}_{\frac{\sqrt{2}c\gamma}{\varphi}X}^{\gamma_\epsilon}[0,S_1]$. {As shown in the proof} of i) in Theorem \ref{TH1},
 $\Theta_2(u)$, $\Sigma_1(u)$, $\Sigma_1(u)$ and $\pk{\sup_{(s,t)\in \eE{D \setminus } E(u)}Z_u(s,t)>m(u)}$ are all negligible in comparison with $J^x(u), x\in (-\IF,\IF]$ and $\pi_1^x(u)$.
 \KD{Hence}
\BQNY
\pk{\sup_{(s,t)\in D_{x,u}}Z_u(s,t)>m(u)}
&\sim&  \sqrt{\pi/a_1}\Phi(x)\mathcal{H}_{\frac{\sqrt{2}c\gamma}{\varphi}X}\mathcal{P}_{\frac{\sqrt{2}c\gamma}{\varphi}X}^{\gamma_\epsilon}\Psi(m(u))\frac{u}
{m(u)\Delta_1(u)}\\
&\sim& \Phi(x)\psi_{\gamma,\IF}(u), \quad u\to \IF.
\EQNY
\KD{In light of (\ref{RT1}), we have} 
 \BQNY
\limit{u} \pk{\frac{\tau_1^*(u)-ut_u}{A(u)}\leq x}=\Phi(x),\ \ x\in(-\IF,\IF].
\EQNY
Further,  (\ref{RT7}) can be proven using the same arguments.
\eE{The joint weak convergence of the passage times follows now by a direct application of} Lemma 2.1 in \cite{HJ13}.\\
{\underline {\it \bf  Case $\eL{\varphi}=\IF$}: The proof of this case follows line by line the same as the proof of case $\varphi=0$ with the exception that we  have to  substitute $B_{\alpha_0}$ with $B_{\alpha_{\IF}}$ throughout the proof of case $\varphi=0$. This completes the proof.\QED

\def\lmu{\CE{\widetilde{m}_u}}


\subsection{ \prooftheo{THF1}}
For any $u$ positive
\BQNY
\psi_{\eL{\gamma}, T}(u)&=&\mathbb{P}\left(\sup_{0\leq t\leq T} \left(X(t)-ct-\gamma\inf_{s\in[0,t]}(X(s)-cs)\right)>u\right)\\
&=&\mathbb{P}\left(\sup_{0\leq t\leq T} Z_{1,u}(s,t)>m_1(u)\right),
\EQNY
where  $m_1(u)=\frac{u+cT}{\sigma(T)}$ and
$$Z_{1,u}(s,t)=\left(\frac{X(t)-\gamma X(s)}{u+c(t-\gamma s)}\right)m_1(u).$$
Below, for notational simplicity we set
$$\sigma_{1,u}^2(s,t):=Var\left(Z_{1,u}(s,t)\right),$$
$$\CE{r_{1}}(s,t,s_1,t_1):=Cor(Z_{1,u}(s,t),Z_{1,u}(s_1,t_1))=Cor(X(t)-\gamma X(s), X(t_1)-\gamma X(s_1)).$$
Let $D_T=\{(s,t), 0\leq s\leq t\leq T\}$ and $A_\delta=[0,\delta]\times [T-\delta , T]$ with $0<\delta<T/2$. For all large $u$
 \BQN\label{F1}
 \pi^*(u)\leq \psi_{\gamma,T}(u)\leq \pi^*(u)+ \pi^{**}(u)
 \EQN
holds  with
 $$\pi^{\eL{*}}(u):=\mathbb{P}\left(\sup_{(s,t)\in A_\delta}Z_{1,u}(s,t)>m_1(u)\right),
\quad \pi^{**}(u):=\mathbb{P}\left(\sup_{(s,t)\in D_T \CE{\setminus}A_\delta}Z_{1,u}(s,t)>m_1(u)\right).$$

 The  idea of the proof is to apply  to $\pi^{\eL{*}}(u)$ Theorem 3.1 in \cite{KEP20151}
 which gives the tail asymptotics of supremum of Gaussian random fields with unique maximum variance point and
to show that $\pi^{**}(u)$ is asymptotically negligible \ehg{compared to $\pi^{*}(u)$.}
For this we need to know the dependence structure of the random field $Z_{1,u}$, which is analyzed in the next step of the proof.

\subsubsection{Dependence structure of $Z_{1,u}$} Proofs of the following lemmas are postponed to Appendix.

\BEL\label{Var} \eL{If}  $\vf$ satisfies {\bf BI} and {\bf BIII}, then the unique maximizer of $\sigma_{1,u}(s,t)$ over
$\{(s,t):0\leq s\leq t\leq T\}$
is $(0,T)$.   Moreover, for $u$ large enough and  as $(s,t)\rw (0,T)$
\BQN\label{Vare}
&&1-\sigma_{1,u}(s,t)\nonumber\\
&& \ \ =\left(\frac{\vfp(T)}{2\vf(T)}-a_3(u)\right)(T-t)(1+o(1))+\left\{
\begin{array}{cc}
\left(\frac{\gamma\vfp(T)}{2\vf(T)}-\gamma a_3(u)\right)s(1+o(1)),& \text{   if  }\vf(s)=o(s)\\
\left(\frac{b(\gamma-\gamma^2)+\gamma\vfp(T)}{2\vf(T)}-\gamma a_3(u)\right)s(1+o(1)),& \text{  if  }\vf(s)\sim b s\\
\frac{\gamma-\gamma^2}{2\vf(T)}\vf(s)(1+o(1)),& \text{  if  }s=o(\vf(s)),
 \end{array}
 \right.
\EQN
where $a_3(u)= \frac{c}{u+cT}\rw 0$, as $u\rw\IF.$
\EEL
\BEL\label{CorF}
If $\vf$ satisfies {\bf BI-BII} and $t^2=o(\vf(t))$ as $t\rw0$, then
\BQN
1-\CE{r_1}(s,t,s_1,t_1)\sim \frac{\vf(|t-t_1|)+\gamma^2\vf(|s-s_1|)}{2\vf(T)}
\EQN
holds for $(s,t), (s_1,t_1)\rw (0,T)$.
\EEL

\subsubsection{Upper estimate of $\pi^{\eL{**}}(u)$}
 By Lemma \ref{Var}, there exists a positive constant $0<\eta<1$ such that
 \BQNY
 \sup_{(s,t)\in \CE{D_T} \setminus A_\delta}Var\left(Z_{1,u}(s,t)\right)\leq  1-\eta.
 \EQNY
 In addition, it follows from {\bf BII} that
 \BQNY
\eL{Var}(Z_{1,u}(s,t)-Z_{1,t}(s',t'))\leq Q_1 \left(|t-t'|^{\alpha_0}+|s-s'|^{\alpha_0}\right), \ \ (s,t)\in D_T.
 \EQNY
\eE{Using \nelem{lemP}  for $u$ large enough we obtain}
 \BQN\label{F2}
 \mathbb{P}\left(\sup_{(s,t)\in \CE{D_T} \setminus A_\delta}Z_{1,u}(s,t)>m_1(u)\right)\leq Q_2T^2(m_1(u))^{\frac{4}{\alpha_0}}\Psi\left(\frac{m_1(u)}{\sqrt{1-\eta }}\right).
 \EQN

\subsubsection{Asymptotics of $\pi^{\eL{*}}(u)$} \mbox{}\\
{\it \underline{Case $s=o(\vf(s))$ as $s\rw 0$}}: In light of Lemma \ref{Var}, \eG{for any positive $\delta$ and $ \epsilon$} sufficiently small we have
 \BQNY
 \mathbb{P}\left(\sup_{(s,t)\in A_\delta}Z_{2,\epsilon}(s,t)>m_1(u)\right)\leq \pi^*(u)\leq \mathbb{P}\left(\sup_{(s,t)\in A_\delta}Z_{2,-\epsilon}(s,t)>m_1(u)\right),
 \EQNY
 where
 $$Z_{2,\pm\epsilon}(s,t)=\frac{\overline{X(t)-\gamma X(s)}}{\left(1+\frac{ \dotvT \pm\epsilon}{2\vf(T)}(T-t)\right)\left(1+\frac{\gamma-\gamma^2\pm\epsilon}{2\vf(T)}\vf(s)\right)},\quad  (s,t)\in A_\delta,$$
\KD{where }$\overline{Z}$ means standardisation of $Z$, i.e., $\overline{Z(t)}= Z(t)/\sqrt{Var(Z(t)}$. In view of Lemma \ref{CorF} and \ehc{using Theorem 3.1} in \cite{KEP20151}, we derive
 \BQN\label{F3}
 \mathbb{P}\left(\sup_{(s,t)\in A_\delta}Z_{2,\pm\epsilon}(s,t)>m_1(u)\right)\sim \mathcal{H}_{B_{\alpha_0}}\mathcal{P}_{B_{\alpha_0}}^{\frac{1-\gamma\pm\epsilon/\gamma}{\gamma}}\frac{2\vf(T)}{ \dotvT \pm\epsilon}
\frac{\Psi\left(\CE{m_1(u)}\right)}{q(u)m_1^2(u)},\ \ u\rw\IF.
 \EQN
 Letting $\delta\rw 0$, $\epsilon\rw0$ leads to
 \BQNY
 \pi^*(u)\sim \mathcal{H}_{B_{\alpha_0}}\mathcal{P}_{B_{\alpha_0}}^{\barg}\frac{2\vf(T)}{ \dotvT }
\frac{\Psi\left(\CE{m_1(u)}\right)}{q(u)m_1^2(u)},\ \ u\rw\IF,
 \EQNY
which together with (\ref{F1}) and (\ref{F2}) establishes the claim.\\
 {\it \underline{Case $\vf(s)\sim bs$ as $s\rw 0$}}:  \ehc{In light of Theorem 3.1 in \cite{KEP20151}}, in this case (\ref{F3}) is changed  \ehb{to}
 \BQNY
 \mathbb{P}\left(\sup_{(s,t)\in A_\delta}Z_{2,\pm\epsilon}(s,t)>m_1(u)\right)\sim \mathcal{P}_{B_{1/2}}^{\frac{ \dotvT \pm \epsilon}{b}}\mathcal{P}_{B_{1/2}}^{\beta(b,\gamma)\pm\frac{ \epsilon}{b\gamma^2}}
\Psi\left(\CE{m_1(u)}\right),\ \ u\rw\IF,
 \EQNY
 with
 $$Z_{2,\pm\epsilon}(s,t)=\frac{\overline{X(t)-\gamma X(s)}}{\left(1+\frac{ \dotvT \pm\epsilon}{2\vf(T)}(T-t)\right)\left(1+\frac{b(\gamma-\gamma^2)+\gamma \dotvT
 \pm\epsilon}{2\vf(T)}s\right)}, \quad  (s,t)\in A_\delta.$$
 Thus letting $\delta\rw 0$, $\epsilon\rw 0$ and using (\ref{F1}) and (\ref{F2}) \ehb{establishes}  the claim.\\
 {\it \underline{ Case $\vf(s)=o(s)$ as $s \rw 0$}}:  For any $\epsilon>0$, if $\delta$ is sufficiently small, then
 \BQNY
 1-\CE{r_1}(s,t,s_1,t_1)\leq \frac{2\left(\vf(|t-t_1|)+\vf(|s-s_1|)\right)}{\vf(T)}\leq \epsilon\left(|t-t_1|+|s-s_1|\right), \ \ (s,t), (s_1,t_1)\in A_\delta.
 \EQNY
 Let $Z^*_{\epsilon}(s,t)$ be a stationary Gaussian field over $[0,T]^2$ with variance $1$ and correlation function
 $$
 \frac{e^{-4\epsilon s}+e^{-4\epsilon t}}{2}, \quad
 s,t\in [0,T].$$
It follows that
 $$1-r_1(s,t,s_1,t_1)<1-\frac{e^{-4\epsilon |s-s_1|}+e^{-4\epsilon |t-t_1|}}{2}, \ \ (s,t), (s_1,t_1)\in A_\delta.$$
 In light of Lemma \ref{Var}, by Slepian's inequality \eG{(see e.g., Theorem 2.2.1 in \cite{AdlerTaylor})} and Theorem 3.1 in \cite{KEP20151}, we have, for \eG{positive} $\delta$ sufficiently small
 \BQN\label{F4}
 \pi^*(u)&\leq& \mathbb{P}\left(\sup_{(s,t)\in A_\delta}\frac{Z^*_{\epsilon}(s,t)}{\left(1+\frac{ \dotvT }{4\CE{\vf(T)}}(T-t)\right)
 \left(1+\frac{\gamma \dotvT }{4\CE{\vf(T)}}s\right)}>m_1(u)\right)\nonumber\\
 &\sim& \mathcal{P}_{B_{1/2}}^{\frac{ \dotvT }{8\epsilon \CE{\vf(T)} }}\mathcal{P}_{B_{1/2}}^{\frac{\gamma \dotvT }{8\epsilon \CE{\vf(T)}}}
\Psi\left(\CE{m_1(u)}\right),\ \ u\rw\IF.
 \EQN
 Moreover,
 \BQNY
  \pi^*(u)\geq \mathbb{P}\left(Z_{1,u}(0,T)>m_1(u)\right)\sim \Psi\left(\CE{m_1(u)}\right),\ \ u\rw\IF.
 \EQNY
 Thus letting $\epsilon\rw 0$ in (\ref{F4}) leads to
 \BQNY
  \pi^*(u)\sim \Psi\left(\CE{m_1(u)}\right),\ \ u\rw\IF,
 \EQNY
 which together with (\ref{F1}) and (\ref{F2}) completes the proof. \QED
 \prooftheo{THF2} For $x>0$, let $T_{x,u}=T-\frac{2\sigma^4(T)x}{ \dotvT u^2}$. For all the three cases, using Theorem \ref{THF1} and Remark \ref{RE} ii) we have
 \BQNY
 \mathbb{P}\left(\frac{ \dotvT u^2(T-\tau_1)}{2\sigma^4(T)}>x\Bigl\lvert \tau_1\leq T\right)=\frac{\psi_{T_{x,u}}(u)}{\psi_T(u)}\sim \frac{\Psi\left(\frac{u+cT_{x,u}}{\sigma(T_{x,u})}\right)}{\Psi\left(\frac{u+cT}{\sigma(T)}\right)}\sim e^{\frac{(u+cT)^2}{2\vf(T)}-\frac{(u+cT_{x,u})^2}{2\vf(T_{x,u})}}, \ \ u\rw\IF,
 \EQNY
 where \ehb{for any $x>0$}
 \BQNY
 \frac{(u+cT)^2}{2\vf(T)}-\frac{(u+cT_{x,u})^2}{2\vf(T_{x,u})}&=& \frac{(u+cT)^2}{2\vf(T)}\left(1-\frac{(1-\frac{c(T-T_{x,u})}{u+cT})^2}{\frac{\vf(T_{x,u})}{\vf(T)}}\right)\\
 &\sim& \frac{(u+cT)^2}{2\vf(T)}\left(1-\frac{(1-\frac{c(T-T_{x,u})}{u+cT})^2}{1-\frac{ \dotvT (T-T_{x,u})}{\vf(T)}}\right)\\
& \rw&
 -x,\quad   u\rw\IF.
 \EQNY
 Thus the claim is established.
 \QED
\section{Appendix}
{
In this \ehb{section we present an extension of Theorem 8.1 in \cite{Pit96}
to threshold-dependent Gaussian fields, followed by a uniform version of Pickands-Piterbarg
lemma motivated by Lemma 2 in \cite{DI2005}.} Then we give lemma that deals with existence and positivity of generalized Piterbarg constants, \rd{which is followed by the proof of (\ref{neg0}).}
\ehb{Finally, we display the proofs of Lemmas 4.1-4.5.} \rd{Before proceeding to the proofs in Appendix, we first present some regularly varying properties on $\sigma^2$. By {\bf AI} and  Theorem 1.7.2 in \cite{BI1989}, we have that
\BQN\label{re1}
\lim_{u\rw\IF}\frac{u\dot{\sigma^2}(u)}{\sigma^2(u)}&= &2\alpha_\IF,\\
  \lim_{u\rw\IF}\frac{u\ddot{\sigma^2}(u)}{\sigma^2(u)}&=&2\alpha_\IF(2\alpha_\IF-1).\label{re2}
\EQN
  Lemma 5.2 in \cite{KrzysPeng2015} shows that {\bf AI} implies that in a neighborhood of $0$
\BQN\label{VarL}
\vf(t)\geq Ct^2,
\EQN
then the function
\BQN\label{A1}
\frac{1}{g_2(t)}=\frac{t^2}{\vf(t)},\quad   t\in(0,\IF)
\EQN
is  regularly varying at infinity with index $2(1-\alpha_\IF)>0$ and is bounded in a neighborhood of zero. By (\ref{A1}) and uniform convergence theorem in \cite{BI1989} we have that for any $T>0$
\BQN\label{re3}\lim_{u\rw\IF}\sup_{t\in (0, T]}\left|\frac{g_2(u)}{g_2(ut)}-|t|^{2-2\alpha_\IF}\right|=0.
\EQN
Moreover,  Potter's bound in \cite{BI1989} shows that for any $0<\epsilon<2\alpha_\IF$, there exists $T>0$ and $Q_1, Q_2>0$ such that for all $s,t>T>0$
\BQN\label{re4}
Q_1\min\left(\left(\frac{t}{s}\right)^{2\alpha_\IF-\epsilon},\left(\frac{t}{s}\right)^{2\alpha_\IF+\epsilon} \right)\leq\frac{\sigma^2(t)}{\sigma^2(s)}\leq Q_2\max\left(\left(\frac{t}{s}\right)^{2\alpha_\IF-\epsilon},\left(\frac{t}{s}\right)^{2\alpha_\IF+\epsilon} \right)
\EQN
By UCT, similarly as in (\ref{re0}) we have that for $0<\lambda<\min(2\alpha_0, 2\alpha_\IF)$ as $u$ sufficiently large,
\BQN\label{re5}
\frac{\sigma^2(ut)}{\sigma^2(u)}=\frac{g_\lambda(ut)}{g_\lambda(u)}|t|^\lambda\leq 2|T|^{2\alpha_\IF-\lambda}|t|^\lambda, \quad t\in[0,T].
\EQN
By {\bf AII} and Theorem 1.7.2 in \cite{BI1989} that
\BQNY
t\dot{\vf}(t)\sim 2\alpha_0\vf(t), \quad t\rw 0,
\EQNY
which combined with (\ref{VarL}) gives that $t/ \dot{\vf}(t)$ is bounded in a neighbourhood of zero.  Therefore if {\bf AI-AII} hold,  we have from (\ref{re1}) that
\BQN\label{A2b}
\ehg{K(t)}:=t ( \dot{\vf}(t))^{-1},\ \ t\in(0,\IF)
\EQN
is a regularly varying function at infinity with index $2(1-\alpha_\IF)>0$ and bounded in a neighbourhood of zero.}

\subsection{Extensions of Piterbarg inequality and Pickands-Piterbarg lemma}
\ehg{Piterbarg inequality, e.g. \cite{Pit96}[Theorem 8.1], provides a precise upper bound for tail distribution of supremum for
a wide class of Gaussian processes.  Our next result extends Piterbarg inequality
to threshold-dependent Gaussian random fields with general covariance structure,
allowing for supremum to be taken on sets that depend on $u$.}

\BEL \label{lemP} Let $X_{u,\tau}(t),t\inr^d, \tau\in K_u, u>0$ be a centered Gaussian field with variance $\sigma_{u,\tau}^2(t), t\in E_{u,\tau} $ and  \he{a.s.}\ continuous sample paths
where $K_u$ are some index sets.  Let further $E_{u,\tau} \subset \prod_{i=1}^d[-M_{u,i},M_{u,i}], u>0,\tau\in K_u$  be compact sets, and put
 $\sigma_{u,\tau}:=\sup_{t\in E_{u,\tau}} \sigma_{u,\tau}(t)$. Suppose that $ 0< a < \sigma_{u,\tau} < b < \IF$ holds for $\tau\in K_u$ and all large $u$.   If   for all $u$ large and for any
$s,t\in E_{u,\tau}$
\BQN\label{sv}
Var\Bigl(X_{u,\tau}(t)- X_{u,\tau}(s) \Bigr)\le 
 C\sum_{i=1}^d\abs{t_i-s_i}^{\eE{\gamma_i}}, \quad s=(s_1 \ldot s_d), \quad t=(t_1 \ldot t_d), \tau\in K_u,
\EQN
with $\ehg{\gamma_i\in (0,2]}, 1\leq i\leq d$,
then for some $C_1>0$ and $u_0>0$ not depending on $u$ and $\tau\in K_u$  
\BQN
\pk{ \sup_{t\in E_{u,\tau}} \abs{X_{u,\tau}(t)}> u}&\le &
C_1 \prod_{i=1}^d\left(M_{u,i}u^{\frac{2}{\gamma_i}}+1\right) \Psi\Bigl(u/\sigma_{u,\tau}\Bigr), \ \ u>u_0.
\EQN
\EEL
\prooflem{lemP} Let $E_{u,\tau}^{(1)}=\{t\in E_{u,\tau}: \sigma_{u,\tau}(t)>\sigma_{u,\tau}/2\}$ and $E_{u,\tau}^{c}:=E_{u,\tau} \setminus E_{u,\tau}^{(1)}$. Then for  $s,t \in E_{u,\tau}^{(1)}$,
$$1-Cor(X_{u,\tau}(t)X_{u,\tau}(s))\leq \frac{Var(X_{u,\tau}(t)- X_{u,\tau}(s))}{2\sigma_{u,\tau}(t)\sigma_{u,\tau}(s)}\leq \frac{2C}{a^2}\sum_{i=1}^d \abs{t_i-s_i}^{\eE{\gamma_i}}.$$
\KD{Let $Y(t), t\in \mathbb{R}^d$} be a homogeneous Gaussian process with variance 1 and correlation function
\eeh{$$r_Y(t)=Cov\left(Y(s), Y(s+t)\right)=e^{-\frac{4C}{a^2}\sum_{i=1}^d \abs{t_i}^{\eE{\gamma_i}}}, \ \ s,t\in \mathbb{R}^d,$$}
\KD{and let} $$L_{\vk{k}}(u)=\prod_{i=1}^d[k_iu^{-\frac{2}{\gamma_i}}, (k_i+1)u^{-\frac{2}{\gamma_i}}],$$
 with
$\vk{k}=(k_1, \dots, k_d)$ and $k_i\in \mathbb{Z}, i=1,\dots, n$.
By \ehb{S}lepian's inequality \eG{(see e.g., Theorem 2.2.1 in \cite{AdlerTaylor})} for $u$ large enough \ehb{we have}
\BQNY
\pk{ \sup_{t\in E_{u,\tau}^{(1)}} \abs{X_{u,\tau}(t)}> u}&\leq& 2\pk{ \sup_{t\in E_{u,\tau}^{(1)}} \overline{X}_{u,\tau}(t)>\frac{ u}{\sigma_{u,\tau}}}\nonumber\\
&\leq& \sum_{\vk{k}:L_{\vk{k}}(u)\cap E_{u,\tau}^{(1)}\neq \emptyset}2\pk{ \sup_{t\in L_{\vk{k}}(u)\cap E_{u,\tau}^{(1)}} \overline{X}_{u,\tau}(t)>\frac{ u}{\sigma_{u,\tau}}}\nonumber\\
&\leq& \sum_{\vk{k}:L_{\vk{k}}(u)\cap E_{u,\tau}^{(1)}\neq \emptyset}2\pk{ \sup_{t\in L_{\vk{k}}(u)\cap E_{u,\tau}^{(1)}} Y(t)>\frac{ u}{\sigma_{u,\tau}}}\nonumber\\
&\leq &\sum_{\vk{k}:L_{\vk{k}}(u)\cap E_{u,\tau}^{(1)}\neq \emptyset}2\pk{ \sup_{t\in L_{\vk{0}}(u)} Y(t)>\frac{ u}{\sigma_{u,\tau}}}.
\EQNY
\eG{Further, by Lemma 6.1 in \cite{Pit96} and the fact that
$$\inf_{\tau\in K_u} \frac{ u}{\sigma_{u,\tau}}\rw \IF, \quad u\rw\IF,$$ we have 
$$\lim_{u\rw\IF}\sup_{\tau\in K_u}\left|\frac{\pk{ \sup_{t\in L_{\vk{0}}(u)} Y(t)>\frac{ u}{\sigma_{u,\tau}}}}{\Psi(\frac{ u}{\sigma_{u,\tau}})}-\prod_{i=1}^d\mathcal{H}_{B_{\gamma_i/2}}[0, \left(4ca^{-2}\right)^{1/\gamma_i}]\right|=0.$$}
Consequently, for $u$ sufficiently large
\BQN\label{FE}
\pk{ \sup_{t\in E_{u,\tau}^{(1)}} \abs{X_{u,\tau}(t)}> u}
&\leq& 2\prod_{i=1}^d\left[\mathcal{H}_{B_{\gamma_i/2}}[0, \left(4ca^{-2}\right)^{1/\gamma_i}]\left([2M_{u,i}u^{\frac{2}{\gamma_i}}]+1\right)\right] \Psi(u/\sigma_{u,\tau})\nonumber\\
&\leq & C_1 \prod_{i=1}^d\left([2M_{u,i}u^{\frac{2}{\gamma_i}}]+1\right) \Psi(u/\sigma_{u,\tau})
\EQN
uniformly with respect to $\tau\in K_u$. \EH{By (\ref{sv}) for any $0<h\leq 1$
 \BQNY
 \sup_{|t_i-s_i|\leq h, 1\leq i\leq d}\sqrt{Var(X_{u,\tau}(t)- X_{u,\tau}(s))}\le 
 \left(C\sum_{i=1}^d\abs{t_i-s_i}^{\eE{\gamma_i}}\right)^{1/2}\leq (Cd)^{1/2} h^{\gamma_0/2},
 \EQNY
 with $\gamma_0=\min_{1\leq i\leq d}\gamma_i$. }
 Thus by  (2.2) in \cite{Breman85} and (\ref{sv}), for any  $E_{u,\tau}^{c}\cap L_{\vk{k}}(1)\neq\emptyset$, we have
\BQN
\pk{ \sup_{t\in E_{u,\tau}^{c}\cap L_{\vk{k}}(1)}|X_{u,\tau}(t)|> \left[b+(2+\sqrt{2})(Cd)^{1/2}\int_1^\IF 2^{-\frac{\gamma_0y^2}{2}}dy\right]x}\leq \frac{5}{2}2^{2d}\sqrt{2\pi}\Psi\left(x\right)
\EQN
for all $x\geq (1+4d\ln 2)^{1/2}$, which implies that  we can find a constant $y$ such that
$$\pk{ \sup_{t\in E_{u,\tau}^{c}\cap L_{\vk{k}}(1)} X_{u,\tau}(t)> y}<1/2$$
holds for all $\vk{k}$ with $E_{u,\tau}^{c}\cap L_{\vk{k}}(1)\neq\emptyset$.
Further, using \ehb{Borell-TIS inequality, see e.g., \cite{AdlerTaylor, Pit20,GennaBorell}}
\BQN\label{FE1}
\pk{ \sup_{t\in E_{u,\tau}^{c}} X_{u,\tau}(t)> u}&\leq& \sum_{\vk{k}: E_{u,\tau}^{c}\cap L_{\vk{k}}(1)\neq \emptyset}\pk{ \sup_{t\in E_{u,\tau}^{c}\cap L_{\vk{k}}(1)} \abs{X_{u,\tau}(t)}> u}\nonumber\\
&\leq &2^d\prod_{i=1}^d(M_{u,i}+1)\Psi\Bigl(2(u-y)/\sigma_{u,\tau}\Bigr)\notag\\
&=&o\left(\prod_{i=1}^d\left(M_{u,i}u^{\frac{2}{\gamma_i}}+1\right) \Psi\bigl(u/\sigma_{u,\tau}\bigr)\right) \notag,
\EQN
hence the claim is established by  considering also
(\ref{FE}).
\QED

\BRM i)  In case $X_{u,\tau}=X$  and $E_{u,\tau} =E$ for all $u$ and $\gamma_i=\gamma, i \le d$, the claim of \nelem{lemP}
coincides with that of Theorem 8.1 in \cite{Pit96}. \ehg{Note in passing that Piterbarg inequality gives sharper bounds than Borell-TIS inequality. The later however holds under weaker assumptions.}\\
\EH{ii) In the formulation of \nelem{lemP} we write $(M_{u,i} u^{2/\gamma_i}+1)$ and not simply
$M_{u,i} u^{2/\gamma_i}$ since we want to cover also the case that $\limit{u} M_{u,i} u^{2/\gamma_i}=0$.}
\ERM

\def\gtu{g_{u,\tau}}
\def\xitu{\xi_{u,\tau}}
\def\chitu{\chi_{u,\tau}}
\def\zetatu{\zeta_{u,\tau}}
\def\vtu{\sigma_{u,\tau}}
\def\bD{E}
\def\sigxiu{\vtu}
\def\cL#1{#1}
\def\H{\mathcal{R}}
\def\coe{C_0(E)}
\ehg{The classical Pickands lemma gives the exact asymptotics of Gaussian processes on short intervals. We present below an extension of that lemma; our result is uniform with respect to some parameter $\tau_u \in K_u$. Let therefore}  $E \subset \R^d$ be a compact set with positive Lebesgue measure containing the origin and let $K_u$ some index sets.
We denote  $\coe$  the space  of all continuous functions $f$ on $E$, \eE{such that $f(0)=0$}, equipped with the sup-norm.
\ehc{For $f_{u,\tau}\in \coe$} define
$$ \xitu(t)=\frac{Z_{u,\tau}(t) }{1+f_{u,\tau}(t)}, \quad t\in E, \ \  \ehc{\tau:=\tau_u} \in K_u,$$
 with $Z_{u,\tau}$ \ehc{a} centered  Gaussian field with unit variance and \he{a.s.}\ continuous sample paths. In the following lemma we derive the uniform asymptotics of
$$ p_{u,\tau}(E):= \pk{\sup_{t\in E}\xi_{u,\tau}(t)> \gtu}, \quad u\to \IF,$$
with respect to $\tau \in K_u$.
We  need the following assumptions, which \ehg{are similar to those imposed in  \cite{KrzysPeng2015}[Lemma 5.1] and \cite{DI2005}[Lemma 2]}.\\
{\bf P1}: $\inf_{\tau\in K_u}g_{u,\tau}\rw\IF$ as $u\rw\IF$.\\
{\bf P2}: Let $\theta_{u,\tau }(s,t)$ be a function such that \BQNY\label{PI2}
\lim_{u\rw\IF}\sup_{\tau \in K_u}\sup_{s\neq t\in E}\left|g_{u,\tau}^2\frac{Var\left(Z_{u,\tau}(t)-Z_{u,\tau}(s)\right)}
{2\theta_{u,\tau }(s,t)}-1\right| \eHHH{=} 0.
\EQNY
There exists \blu{a centered Gaussian random field $V(t),t\inr^d$}  with $V(0)=0$,
covariance function $(\sigma_V^2(t)+ \sigma_V^2(s)-\sigma_V^2(t-s))/2,s,t\inr^d$ and \he{a.s.}\ continuous sample paths such that
$$\limit{u}\sup_{\tau \in K_u}|\theta_{u, \tau}(s, t)-\sigma_V^2(t-s)|= 0, \quad \forall s, t\in E.$$\\
{\bf P3}: There exists some $\ehb{a} >0$ such that
\BQNY
\limsup_{u\rw\IF}\sup_{\tau \in K_u}\sup_{s\neq t, s,t\in E}\frac{\theta_{u, \tau }(s,t)}{\sum_{i=1}^{d}|s_i-t_i|^{\ehb{a}}}<\IF
\EQNY
and further
\BQNY
\lim_{\epsilon\downarrow  0}\limsup_{u\rw\IF}\sup_{\tau \in K_u}\sup_{\norm{t-s}<\epsilon,s,t\in E}g_{u,\tau}^2
\EE{ \left[ Z_{u,\tau}(t)-Z_{u,\tau }(s)\right]Z_{u,\tau}(0)}=0.
\EQNY

\BEL \label{PICKANDS}
 Let $g_{u,\tau}, V, \theta_{u,\tau}$ satisfy {\bf P1-P3}.
  Assume that $f_{u,\tau}\in  \coe ,u>0,\tau \in K_u$
\BQN\label{PI1}
 \lim_{u\to \IF} \sup_{\tau\in K_u, t\in E} \abs{ g_{u,\tau}^2 f_{u,\tau}(t)- f(t)}=0.
 \EQN
Then we have
\BQN
\ehc{\limit{u} \sup_{\tau \in K_u} }\ABs{
 \frac{p_{u,\tau}(E)}{\Psi(g_{u,\tau})} -  \H_{V}^{f}(\bD)}=0,
\EQN
with
$\H_{V}^{f}(\bD):=\EE{\sup_{t\in E}e^{\sqrt{2}V({t})-\sigma_V^2({t})-f({t})}}\in (0,\IF)$.
\EEL

\prooflem{PICKANDS} By conditioning on $\xitu(0)=g_{u,\tau}-\frac{w}{g_{u,\tau}}, w\in\mathbb{R}$
for all $u>0$ large we obtain
\BQNY
 { \sqrt{2 \pi} g_{u,\tau} e^{\frac{g_{u,\tau}^2}2} }  \pk{\sup_{t\in E}\xitu(t)>\gtu} =
\int_{\R} e^{w-\frac{w^2}{2\gtu^2}}\pk{\sup_{t\in E}\chitu(t)> w}\, dw=:
I_{u,\tau},
\EQNY
where
$$\chitu(t)=\zetatu (t)|\zetatu({0})=0, \quad \zetatu(t)= \gtu(\xitu({t})-\gtu)+w.$$
In order to establish the proof we need to show that
\BQN \label{esB}
\limit{u} \sup_{\tau \in K_u} \ABs{I_{u,\tau}- \H_{V}^{f}(\bD)} = 0.
\EQN
\ehb{It follows} that
\BQNY
 \sup_{\tau \in K_u} \abs{I_{u,\tau}- \H_{V}^{f}(\bD)}&\leq&  \sup_{\tau \in K_u} \ABs{\int_{-M}^M
\Biggl[ e^{w-\frac{w^2}{2 \gtu^2}}\pk{\sup_{t\in E}\chitu(t)> w}- e^w\pk{\sup_{t\in E}V(t)> w}\Biggr]\, dw }\\
&& + \sup_{\tau \in K_u}\lefteqn{\int_{\abs{w}> M } e^{w-\frac{w^2}{2 \gtu^2}} \pk{\sup_{t\in E}\chitu(t)> w}\, dw}\\
&& +\int_{\abs{w}> M} e^{w}\pk{\sup_{t\in E}V(t)> w}dw.
\EQNY
Next, we give \kk{an} upper bound of each term in the right hand side of the above inequality.
Clearly, $\chitu(0)=0$ almost surely, and the finite-dimensional distributions of
$\chitu(t), t\in \bD$ coincide with that of
\BQNY
\frac{1}{1+f_{u,\tau}(t)}\left(\gtu\Bigl(Z_{u,\tau}({t})-r_{Z_{u,\tau}}({t},{0})Z_{u,\tau}({0})\Bigr)+\mu_{u,\tau,w}(t)\right), \quad {t}\in\bD\ehb{,}
\EQNY
where
$$\mu_{u,\tau,w}(t)=-g_{u,\tau}^2\left(1-r_{Z_{u,\tau}}(t,0)+f_{u,\tau}(t)\right)+w(1-r_{Z_{u,\tau}}(t,0)+f_{u,\tau}(t)), \ \ r_{Z_{u,\tau}}(t,s):=Cor(Z_{u,\tau}(t),Z_{u,\tau}(s)).$$

Consequently, by {\bf P1-P3} \KD{and} (\ref{PI1}) {we have that} uniformly with respect to $t\in\bD, \tau \in K_u, w\in [-M,M]$
\BQN\label{eqex}
\mu_{u,\tau,w}(t)
& \to & -(\sigma_V^2({t})+f({t})), \  u\rw\IF
\EQN
and also for any $(s,t) \in \bD$ uniformly with respect to $ \tau \in K_u, w\in [-M,M]$
\BQN\label{eqco}
v_u(s,t)&:=&Var\Bigl((1+f_{u,\tau}(t))\chitu({t})-(1+f_{u,\tau}(s))\chitu({s})\Bigr)\nonumber\\
&=&\gtu^2\LT[\EE{\Bigl( Z_{u,\tau}({t})-Z_{u,\tau}({s})\Bigr)^2}-\LT(r_{Z_{u,\tau}}({t},{0})-r_{Z_{u,\tau}}({s},{0})\RT)^2\RT]\notag \nonumber\\
& \to & 2 Var(V({t})-V({s})), \  u\rw\IF.
\EQN
Note that $v_u(s,t) $ does not depend on $w$ and $f\in \coe$. Consequently,
following the proof of Lemma 4.1 in \cite{Yimin15},
{the} finite-dimensional distributions of $(1+f_{u,\tau}(t))\chitu({t})$ converge uniformly for $\tau \in K_u, w\in [- M,M]$ where $M>0$ is fixed.
By {\bf P3}, the uniform convergence in \eqref{eqex}, (\ref{eqco}) and
\BQN\label{fu}\lim_{u\rw \IF}\sup_{\tau \in K_u, t\in \bD}|f_{u,\tau}(t)|=0
\EQN
imply along the lines of the proof of second part of Lemma 4.1 in \cite{Yimin15} that for arbitrary $M>0, \ve>0$
$$\limit{u}\sup_{\tau \in K_u, w\in [-M,M], \EH{w\not \in [-\ve, \ve]}}\ABs{\pk{\sup_{t\in E}\chitu(t)> w}- \pk{\sup_{t\in E}V(t)> w} }=0  ,   $$
\EH{where we use the fact that $\sup_{t\in E} V(t)$ has a continuous distribution
$H(t),t\ge 0$ for all $t>0$, see e.g., \cite{AZI}[Theorem 7.1] (recall that since $0\in E$ and $V(0)=0$, then $\sup_{t\in E} V(t)\ge 0$). Further, by ${\bf P1}$}
$$\limit{u}\sup_{\tau \in K_u, w\in [-M,M]}e^w [1-   e^{-\frac{w^2}{2 \gtu^2}}] \le
 \frac{e^M M^2}{2 \liminf_{u\to \IF }\inf_{\tau \in K_u}  \gtu^2} \to 0, \quad u \to \IF      $$
\EH{we obtain by the fact that $\ve>0$ can be chosen arbitrary small}
$$\limit{u}\sup_{\tau \in K_u}\ABs{\int_{-M}^M
\Biggl[ e^{w-\frac{w^2}{2 \gtu^2}}\pk{\sup_{t\in E}\chitu(t)> w}- e^w\pk{\sup_{t\in E} V(t)> w}\Biggr]\, dw }=0.     $$

Using \eqref{eqex} for $\delta \in (0,1/2)$, $|w|>M$ with $M$ sufficiently large and all $u$ large we have
$$ \sup_{\tau \in K_u,t\in E} \EE{(1+f_{u,\tau}(t))\chitu({t})} \le \delta \abs{w}.$$
It follows from {\bf P3} that for $u$ large enough,
$$ Var\left((1+f_{u,\tau}(t))\chitu({t})-(1+f_{u,\tau}(s))\chitu({s})\right)\leq Q\sum_{i=1}^{d}|s_i-t_i|^{\ehb{a}}, \quad (s,t)\in E.$$
Thus by (\ref{fu}) and the result of \nelem{lemP}, we obtain for some $\ve,\delta \in (0,1/2)$ and all $u$ and $M$ large
\BQNY
\lefteqn{\int_{\abs{w}> M } e^{w-\frac{w^2}{2 \gtu^2}} \pk{\sup_{t\in E}\chitu(t)> w}\, dw}\\
&\le &
\int_{\abs{w}> M } e^{w}\pk{\sup_{t\in E}  \left((1+f_{u,\tau}(t))\chitu({t})- \EE{(1+f_{u,\tau}(t))\chitu({t})}\right)> w/2-\sup_{t\in E, \tau \in K_u} \EE{(1+f_{u,\tau}(t))\chitu({t}) }}\, dw\\
&\leq& \int_{-\IF}^{-M} e^{w}dw+\int_M^\IF e^w \pk{\sup_{t\in E}  (\chitu(t)- \EE{\chitu(t)})> w/2-\delta w }\, dw\\
&\le & e^{-M}+ \int_{M}^\IF  e^{w}\Psi \bigl( (1- \ve)(1/2-\delta)w \bigr) \, dw \\
&=:&A_1(M) \to   0, \quad M\to  \IF.
\EQNY
Moreover, using \ehg{Borell-TIS} inequality \kk{(see e.g., \cite{AdlerTaylor, GennaBorell})}
$$ A_2(M):= \int_{\abs{w}> M} e^{w}\pk{\sup_{t\in E} V(t)> w}\, dw \to 0, \quad M \to \IF.$$
Hence \eqref{esB} follows since
\BQNY
 \sup_{\tau \in K_u} \abs{I_{u,\tau}- \H_{V}^{f}(\bD)}&\le &
 \sup_{\tau \in K_u} \ABs{\int_{-M}^M
\Biggl[ e^{w-\frac{w^2}{2\gtu^2}}\pk{\sup_{t\in E}\chitu(t)> w}- e^w\pk{\sup_{t\in E}V(t)> w}\Biggr]\, dw }
+ A_1(M)+ A_2(M)\\
&\to & A_1(M)+ A_2(M), \quad u\to \IF, \\
&\to &0, \quad M\to \IF.
\EQNY
Since further
$$\lim_{u\rw\IF}\sup_{\tau\in K_u}\left|\sqrt{2 \pi} g_{u,\tau} e^{\frac{g_{u,\tau}^2}2}\Psi(g_{u,\tau})-1 \right|=0,$$
the proof \kk{is completed.} \QED

\subsection{Piterbarg-type constant}
In this subsection we prove the existence and positivity of the generalized Piterbarg constant that appears in Theorem \ref{TH1}.
Let $X$ be a centered Gaussian process with stationary increments, \he{a.s.}\ continuous sample paths \KD{and variance function satisfying} the following two assumptions:\\
{\bf C0:} $\sigma^2(t)$ is regularly varying at infinity with index $2\alpha_\IF\in (0,2)$ and \ehg{its first derivative is} continuously differentiable over $(0,\IF)$ with $\dot{\sigma^2}(t)$ being ultimately monotone at infinity.\\
{\bf C1:} $\sigma^2(t)$ is regularly varying at zero with index $2\alpha_0\in (0,2]$.\\
Then we have
\BQNY
1-Cor\left(X(ut), X(us)\right)=\frac{\sigma^2(u|t-s|)-(\sigma(ut)-\sigma(us))^2}{2\sigma(ut)\sigma(us)}
=\frac{\sigma^2(u|t-s|)-(u\dot{\sigma}(u\theta)(t-s))^2}{2\sigma(ut)\sigma(us)},
\EQNY
with $\theta\in [s,t]$. Note that (\ref{re1}) implies
$$\lim_{u\rw\IF}\frac{u\dot{\sigma}(u)}{\sigma(u)}=\alpha_\IF,$$
 which together with (\ref{re3}) implies that
\BQN\label{New}
1-Cor\left(X(ut), X(us)\right)&\sim&\frac{\sigma^2(u|t-s|)}{2\sigma(ut)\sigma(us)}\left(1-\frac{\alpha_\IF^2}{\theta^2}\frac{\sigma^2(u\theta)
(t-s)^2}{\sigma^2(u|t-s|)}\right)\nonumber\\
&=&\frac{\sigma^2(u|t-s|)}{2\sigma(ut)\sigma(us)}\left(1-\alpha_\IF^2\frac{g_2(u\theta)}{g_2(u|t-s|)}\right)\\
&\sim &\frac{\sigma^2(u|t-s|)}{2\sigma^2(u)}
\EQN
 as $u\rw\IF$ for $s,t\in [1,1+u^{-1}\ln u]$.
\BEL\label{EXISTENCE}
\ehc{If  $X$ is a centered} Gaussian process with stationary increments \ehc{ and \he{a.s.}\ continuous sample paths such that its}
 variance function satisfies {\bf C0},{\bf C1},  then
$$\mathcal{P}_{X}^{ a}=\lim_{S\rw\IF}\mathcal{P}_{ X}^{a }[0,S]<\IF$$
holds for any $a \in (0,\IF)$.
\EEL
\prooflem{EXISTENCE} \EH{We first introduce some notation. For $S>0, u>1$ define}
$$Y_u(t)=\frac{\overline{X}(u(t+1))}{1+\frac{a\sigma^2(ut)}{2\sigma^2(u)}}, \quad t\in [0, u^{-1}\ln u], $$
$$I_k(u)=[ku^{-1}S, u^{-1}(k+1)S ], \quad 0\leq k\leq N(u), \EH{\text{  with } N(u):=}[S^{-1}\ln u]+1.$$ It follows that for $S$ sufficiently large
\BQN\label{New1}
p_0(u)\leq \pk{\sup_{t\in [0, u^{-1}\ln u]}Y_u(t)>\sqrt{2}\sigma(u)}\leq p_0(u)+\sum_{k=1}^{N(u)}p_k(u),
\EQN
 where
$$p_0(u)=\pk{\sup_{t\in I_0(u)}Y_u(t)>\sqrt{2}\sigma(u)}, $$
$$ p_k(u)=\pk{\sup_{t\in I_k(u)}\overline{X}(u(t+1))>\sqrt{2}\sigma(u)\left(1+\frac{a\sigma^2(kS)}{4\sigma^2(u)}\right)}, \quad  k\geq 1. $$
In order to apply Lemma \ref{PICKANDS}, by (\ref{New}) we set
\BQN K_u=\{k: 0\leq k\leq N(u)\}, \ \ E=[0,S], \ \ g_{u,k}=\sqrt{2}\sigma(u)\left(1+\frac{a\sigma^2(kS)}{4\sigma^2(u)}\right), k\in K_u,
\EQN $$Z_{u,k}(t)=\overline{X}(u(u^{-1}kS+u^{-1}t+1)), \quad k\in K_u,$$
$$\theta_{u,k}(s,t)=g_{u,k}^2\frac{\sigma^2(|t-s|)}{2\sigma^2(u)},\quad  s,t \in E, k\in K_u,$$
$$f_{u,0}(t)=\frac{a\sigma^2(t)}{2\sigma^2(u)},\quad  t\in E, \ \ f_{u,k}=0, \quad k\in K_u \setminus \{0\}, \ \ V=X. $$
\ehd{Since {\bf P1-P2} are obviously fulfilled, we shall verify} next {\bf P3}. By {\bf C1} we have, for $u$ sufficiently large
$$\theta_{u,k}(s,t)=g_{u,k}^2\frac{\sigma^2(|t-s|)}{2\sigma^2(u)}\leq 2\sigma^2(|t-s|)\leq Q|t-s|^{\alpha_0}, \ \ s,t\in E, k\in K_u.$$
Moreover, by (\ref{New})
\BQNY
&&\sup_{k \in K_u}\sup_{\norm{t-s}<\epsilon,s,t\in E}g_{u,k}^2
\EE{ \left[ Z_{u,k}(t)-Z_{u,\tau }(s)\right]Z_{u,k}(0)}\\
&& \ \ \leq \sup_{k \in K_u}\sup_{\norm{t-s}<\epsilon,s,t\in E}g_{u,k}^2
\left(\frac{\sigma^2(t)}{2\sigma^2(u)}(1+o(1))-\frac{\sigma^2(s)}{2\sigma^2(u)}(1+o(1))\right)\\
&& \ \ \leq \sup_{k \in K_u}\sup_{\norm{t-s}<\epsilon,s,t\in E}\frac{g_{u,k}^2}{2\sigma^2(u)}
\left(\abs{\sigma^2(t)-\sigma^2(s)}+o(1)\right)\rw 0, \ \  u\rw\IF, \epsilon\downarrow 0.
\EQNY
Thus {\bf P3} is satisfied. \kk{Hence}
$$ g_{u,0}^2f_{u,0}(t)\rw a\sigma^2(t), \quad u\rw\IF
$$
uniformly with respect to $t\in E$ and
$$g_{u,k}^2f_{u,k}(t)=0, \ \ t\in E, k\in K_u \ehd{\setminus} \{0\},  \quad u>0,
$$
implying  that
\BQNY
\limit{u} \frac{p_0(u)}{\Psi(\sqrt{2}\sigma(u))}=  \mathcal{P}_X^a[0,S]
\EQNY
and
\BQN
\lim_{u\rw\IF}\sup_{k\in  K_u/\{0\}}\ABs{\frac{p_k(u)}{\Psi\left(\sqrt{2}\sigma(u)\left(1+\frac{a\sigma^2(kS)}{4\sigma^2(u)}\right)\right)}-\mathcal{H}_X[0,S]}=0.
\EQN
Dividing (\ref{New1}) by $\Psi(\sqrt{2}\sigma(u))$, letting $u\rw\IF$  and applying (\ref{re4}) \EH{ for  sufficiently large $S_1$}
we have
\BQNY
\mathcal{P}_X^a[0,S]&\leq& \mathcal{P}_X^a[0,S_1]+\mathcal{H}_X[0,S_1]\sum_{k=1}^{\IF} e^{-\frac{a\sigma^2(kS_1)}{2}}\\
&\leq& \mathcal{P}_X^a[0,S_1]+\mathcal{H}_X[0,S_1]\sum_{k=1}^{\IF} e^{-Q_1(kS_1)^{\alpha_\IF}}\\
&\leq & \mathcal{P}_X^a[0,S_1]+\mathcal{H}_X[0,S_1] e^{-Q_2S_1^{\alpha_\IF}}.
\EQNY
Letting $S\rw\IF$ leads to
\BQNY
\lim_{S\rw\IF}\mathcal{P}_X^a[0,S]\leq \mathcal{P}_X^a[0,S_1]+\mathcal{H}_X[0,S_1] e^{-Q_2S_1^{\alpha_\IF}}<\IF
\EQNY
establishing the proof.\QED
\subsection{Proofs of (\ref{neg0})}
We begin with $p_1$, assuming that $T\in \mathbb{N}$ is sufficiently large.  For $(s,t)\in [k,k+1]\times[l,l+1]$ with $l\geq T$ and $0\leq k\leq l$, \rd{ by (\ref{re4})}, we have
\BQNY
Var(Z_u(s,t))&=&\frac{[(1-\gamma)\sigma^2(ut)+(\gamma^2-\gamma)\sigma^2(us)+\gamma\sigma^2(u|t-s|)](1+ct_u)^2}{\left(1+c(t-\gamma s)\right)^2\sigma^2(ut_u)}\\
&\leq&Q\frac{t^{2\alpha_\IF+\epsilon}}{(1+c(1-\gamma)t)^2}\\
&\leq &Qt^{-(2-2\alpha_\IF-\epsilon)}\\
&\leq &
 Ql^{-(2-2\alpha_\IF-\epsilon)}
\EQNY
 for $u$ sufficiently large, with $0<\epsilon<\min(2\alpha_\IF, 2-2\alpha_\IF)$.
Moreover, in view of (\ref{re5}), for $(s,t), (s_1,t_1)\in [k,k+1]\times[l,l+1]$ with $l\geq T$ and $0\leq k\leq l$ and $u$ large enough
\BQNY
Var(\overline{Z}_u(s,t)-\overline{Z}_u(s_1,t_1))&=&2-2Cov\left(\frac{X(ut)-\gamma X(us)}{\sigma_\gamma(us,ut)},\frac{X(ut_1),\gamma X(us_1)}{\sigma_\gamma(us_1,ut_1)}\right)\\
&=&\frac{Var\left(X(ut)-X(ut_1)+\gamma X(us_1)-\gamma X(us)\right)-(\sigma_\gamma(us,ut)-\sigma_\gamma(us_1,ut_1))^2}{\sigma_\gamma(us,ut)\sigma_\gamma(us_1,ut_1)}\\
&\leq& 2\frac{\sigma^2(u|t-t_1|)+\sigma^2(u|s-s_1|)}{\sigma_\gamma(us,ut)\sigma_\gamma(us_1,ut_1)}\\
&\leq&  \frac{4}{(1-\gamma)^2}\frac{\sigma^2(u|t-t_1|)+\sigma^2(u|s-s_1|)}{\sigma^2(ul)}\\
&\leq&Q_T\left(|s-s_1|^{ \CE{\lambda} }+|t-t_1|^{ \CE{\lambda} }\right),
\EQNY
where $Q_T$ is a positive constant depending on $T$ and $0< \CE{\lambda} <\min(2\alpha_0, 2\alpha_\IF)$.  Thus from the above results
\eE{and using further \nelem{lemP}},
for $T$ large enough we have
\BQNY
p_1(u)&\leq& \sum_{l=T}^\IF\sum_{k=0}^l\pk{\sup_{(s,t)\in [k,k+1]\times[l,l+1]}Z_u(s,t)>m(u)}\\
&\leq& \sum_{l=T}^\IF\sum_{k=0}^l\pk{\sup_{(s,t)\in [0,1]^2}\overline{Z}_u(s+k,t+l)>\frac{m(u)}{\sqrt{Ql^{-(2-2\alpha_\IF-\epsilon)}}}}\\
&\leq& \sum_{l=T}^\IF Q_2l \left(m^2(u)l^{2-2\alpha_\IF-\epsilon}\right)^{2/\lambda} e^{-Q_1m^2(u)l^{2-2\alpha_\IF-\epsilon}}\\
&\leq& e^{-Q_3m^2(u)T^{2-2\alpha_\IF-\epsilon}}\\
&=&o\left(\frac{u}{m(u)\Delta_1(u)}\Psi(m(u))\right).
\EQNY
Next, we show that $p_2(u)$ is also negligible. By UCT, we have
\BQNY
Var(Z_u(s,t))\rw \frac{[(1-\gamma)t^{2\alpha_\IF}+(\gamma^2-\gamma)s^{2\alpha_\IF}+\gamma|t-s|^{2\alpha_\IF}](1+c\eL{t_*})^2}{\left(1+c(t-\gamma s)\right)^2\eL{t_*}^{2\alpha_\IF}}=\frac{f(s,t)}{f(0,\eL{t_*})}, \ \ u\rw\IF
\EQNY
 uniformly over $D_{\delta,u}$, where $f(s,t)$ is defined in (\ref{eq1}) with $(0,\eL{t_*})$ the unique maximum point over $D$. Consequently, there exists a constant $0<b_\delta<1$ such that for $u$ large enough
\BQNY
\sup_{(s,t)\in D_{\delta, u}}Var(Z_u(s,t))<\ehg{b_\delta}.
\EQNY
Furthermore, by (\ref{re5}) for $u$ large enough we have
\BQN\label{NEG}
Var(Z_u(s,t)-Z_u(s_1,t_1))&=&\frac{(1+ct_u)^2}{\sigma^2(ut_u)}\ehg{\mathbb{E}\Bigl\{} \Bigl(\frac{X(ut)-\gamma X(us)}{1+c(t-\gamma s)}-\frac{X(ut_1)-\gamma X(us_1)}{1+c(t_1-\gamma s_1)}\Bigr)^2\ehg{\Bigr\}}\nonumber\\
&\leq& Q_4\left(\frac{\sigma^2(u|t-t_1|)}{\sigma^2(ut_u)}+
\frac{\sigma^2(u|s-s_1|)}{\sigma^2(ut_u)}+(t-t_1)^2+(s-s_1)^2\right)\nonumber\\
&\leq&Q_5(|t-t_1|^{ \CE{\lambda} }+|s-s_1|^{ \CE{\lambda} }), \ \ (s,t), (s_1,t_1)\in D_T,
\EQN
with $\lambda\in (0,\min (2\alpha_0, 2\alpha_\IF))$. Consequently, by \nelem{lemP}
$$p_2(u)\leq Q_6T^2(m(u))^{4/\lambda} \Psi\left(\frac{m(u)}{b_\delta}\right)=o\left(\frac{u}{m(u)\Delta_1(u)}\Psi(m(u))\right).$$
Finally, we focus on $p_3(u)$. In light of Lemma \ref{L1}, we know that for $\delta$ sufficiently small and $u$ sufficiently large
\BQNY
\sup_{(s,t)\in D_{\delta,u}^*}Var(Z_u(s,t)) &\leq& \sup_{(s,t)\in D_{\delta,u}^*}\left(1-\frac{a_1}{2}(t-t_u)^2-\frac{a_2}{2}\frac{\sigma^2(us)}{\sigma^2(u)}\right)\\
&\leq &\sup_{(s,t)\in D_{\delta,u}^*}\left(1-\frac{a_1}{2}(t-t_u)^2\right)\\
&\leq &1-Q_7\left(\frac{\ln m(u)}{m(u)}\right)^2,
\EQNY
which together with (\ref{NEG}) and the application of \nelem{lemP}  leads to
$$p_3(u) \leq Q_8(m(u))^{\frac{4}{ \CE{\lambda} }}\Psi\left(\frac{m(u)}{\sqrt{1-Q_7\left(\frac{\ln m(u)}{m(u)}\right)^2}}\right)=o\left(\frac{u}{m(u)\Delta_1(u)}\Psi(m(u))\right), \quad u\to \IF$$
\EH{establishing}  (\ref{neg0}). \QED

\subsection{Proofs of Lemmas \ref{L1}-\ref{CorF}}
In this section we present details of the proof of Lemmas \ref{L1}-\ref{CorF}.\\
 \prooflem{L1} For any $u>0$ we have
\BQNY
\sigma^2_{\gamma,u}(s,t)=\frac{(1-\gamma)\sigma^2(ut)+(\gamma^2-\gamma)\sigma^2(us)+\gamma\sigma^2(u(t-s))}{(1+c(t-\gamma s))^2}.
\EQNY
By UCT we have as $u\to \IF$
\BQN\label{eq1}
\frac{\sigma^2_{\gamma,u}(s,t)}{\sigma^2(u)}\rw\frac{(1-\gamma)t^{2\alpha_\IF}+(\gamma^2-\gamma)s^{2\alpha_\IF}+\gamma (t-s)^{2\alpha_\IF}}{(1+c(t-\gamma s))^2}=:f(s,t)
\EQN
uniformly for $0\leq s\leq t\leq T$ with $T$ any positive constant. \rd{Using (\ref{re4})}  for any $0<\epsilon<\min(2\alpha_\IF,2-{2\alpha_\IF})$ there exists a constant $u_\epsilon>0$ such that for all $0\leq s\leq t<\IF, t>T>1$  and $u>u_\epsilon$, we have
\BQN\label{eq2}
\frac{\sigma^2_{\gamma,u}(s,t)}{\sigma^2(u)}\leq Q\frac{((1-\gamma)t^{2\alpha_\IF+\epsilon}+\gamma t^{2\alpha_\IF+\epsilon})}{(1+c(t-\gamma s))^2}\leq \frac{Q}{t^{2-2\alpha_\IF-\epsilon}}\RW 0, \ \ t\rw\IF.
\EQN
\rd{It follows from \cite{HA2013} that $f(s,t)$ has one unique maximum point $(0,\eL{t_*})$ over $D$, which combined with (\ref{eq1}) and (\ref{eq2}) yields that }for $u$ large enough, the maximum point of $\sigma_{\gamma,u}^2(s,t)$ denoted by $(s_u,t_u)$ must be attained over $0\leq s\leq t\leq T$ with $T>\eL{t_*}$ large enough. Further, $(s_u,t_u)\rw(0,\eL{t_*})$. By contradiction, suppose that $(s_u,t_u)\rw (s_1^*,t_1^*)\neq (0,\eL{t_*})$. Hence, by (\ref{eq1}), we have that
\BQNY
f(s_1^*,t_1^*)=\lim_{u\rw\IF}\frac{\sigma^2_{\gamma,u}(s_u,t_u)}{\sigma^2(u)}\geq \lim_{u\rw\IF}\frac{\sigma^2_{\gamma,u}(0,\eL{t_*})}{\sigma^2(u)}=f(0,\eL{t_*}).
\EQNY
This contradicts the fact that  $(0,\eL{t_*})$ is the unique maximum point of $f(s,t)$ over $D$.
  Next,  we prove that the maximum point is unique. It follows that for $0<s<t<\IF$
\BQN\label{eq:derivative}
\frac{\partial \sigma_{\gamma,u}^2(s,t)}{\partial s}
&=&\gamma A^{-4}(s,t)\left\{\left((\gamma-1)\dot{\sigma^2}(us)u-\dot{\sigma^2}(u(t-s))u\right)A^2(s,t)+2c\sigma_\gamma^2(us,ut)A(s,t)\right\},\nonumber\\
\frac{\partial \sigma_{\gamma,u}^2(s,t)}{\partial t}
&=&A^{-4}(s,t)\left\{\left((1-\gamma)\dot{\sigma^2}(ut)u+\gamma\dot{\sigma^2}(u(t-s))u\right)A^2(s,t)-2c\sigma_\gamma^2(us,ut)A(s,t)\right\},
\EQN
with $A(s,t)=1+c(t-\gamma s)$.
Suppose that $s_u>0$, then by the continuous differentiability of $\sigma^2_{\gamma, u}(s,t)$, we have
$$\frac{\partial \sigma_{\gamma,u}^2(s,t)}{\partial s}\Bigl\lvert_{(s,t)=(s_u,t_u)}=\frac{\partial \sigma_{\gamma,u}^2(s,t)}{\partial t}\Bigl\lvert_{(s,t)=(s_u,t_u)}=0,$$
which implies that
\BQNY
\dot{\sigma^2}(us_u)=\dot{\sigma^2}(ut_u)-\dot{\sigma^2}(u(t_u-s_u))=\ddot{\sigma^2}(u\theta_u)us_u,
\EQNY
with $\theta_u\in(t_u-s_u,t_u)$.
\ehg{For $K(t)= t / \dot{\vf}(t)$ defined in (\ref{A2b})}, the last equation can be re-written as
\BQN\label{eq3}
\frac{u\theta_u\ddot{\sigma^2}(u\theta_u)}{\dot{\sigma^2}(u\theta_u)}\frac{K(us_u)}{K(u\theta_u)}=1.
\EQN
Since {\bf AI} holds, then by (\ref{re1}-\ref{re2}) and using UCT, we have
\BQNY
\ehc{\limit{u}}\frac{u\theta_u\ddot{\sigma^2}(u\theta_u)}{\dot{\sigma^2}(u\theta_u)}= 2\alpha_\IF-1,\quad
\ehc{\limit{u}}\frac{K(us_u)}{K(u\theta_u)}= 0.
\EQNY
Hence, for $u$ large enough
\BQNY
\frac{u\theta_u\ddot{\sigma^2}(u\theta_u)}{\dot{\sigma^2}(u\theta_u)}\frac{K(us_u)}{K(u\theta_u)}<1,
\EQNY
which contradicts  (\ref{eq3}). Consequently, for $u$ large enough then  $s_u=0$ and $t_u$ is the maximum point of $\sigma_{\gamma,u}^2(0,t)=\frac{\sigma^2(ut)}{(1+ct)^2}$. It follows that (the following derivatives are all with respect to $t$)
\BQNY
\frac{\dot{\sigma_{\gamma,u}^2}(0,t)}{\sigma^2(u)}\rw \dot{f}(0,t), \ \ \mbox{and} \ \ \frac{\ddot{\sigma_{\gamma,u}^2}(0,t)}{\sigma^2(u)}\rw \ddot{f}(0,t)<0, \ \ u\rw\IF
\EQNY
hold uniformly over $[\eL{t_*}-\delta,\eL{t_*}+\delta]$ for $\delta>0$ small enough. This implies that $\frac{\dot{\sigma_{\gamma,u}^2}(0,t)}{\sigma^2(u)}$ is decreasing over $[\eL{t_*}-\delta,\eL{t_*}+\delta]$ for $\delta>0$. Thus $t_u$ is unique and then $(0,t_u)$ is unique.
We also have that the first derivative of $\sigma_{\gamma, u}^2(0,t)$ with respect to $t$ at point $t_u$ equals zero (see (\ref{eq:derivative})), i.e.,
 $$\frac{\partial \sigma_{\gamma,u}^2(0,t)}{\partial t}\Big|_{t=t_u}
=A^{-4}(0,t_u)\left\{\left((1-\gamma)\dot{\sigma^2}(ut_u)u+\gamma\dot{\sigma^2}(ut_u)u\right)A^2(0,t_u)-2c\sigma_\gamma^2(0,ut_u)A(0,t_u)\right\}=0,$$
 which is equivalent to
\BQN\label{eq4}
u\dot{\sigma^2}(ut_u)(1+ct_u)^2=2c\sigma^2(ut_u)(1+ct_u).
\EQN
For any $u>0$
\BQNY
&&(1+ct_u+c(t-t_u-\gamma s))^2\sigma^2(ut_u)\nonumber\\
&\quad&\quad=(1+ct_u)^2\sigma^2(ut_u)+2c(1+ct_u)(t-t_u-\gamma s)\sigma^2(ut_u)+c^2(t-t_u-\gamma s)^2\sigma^2(ut_u)
\EQNY
\ehg{and by Taylor expansion }
\BQNY \sigma^2(ut)=\sigma^2(ut_u)+\dot{\sigma^2}(ut_u)u(t-t_u)+\frac{1}{2}\ddot{\sigma^2}(u\theta_{1,u})u^2(t-t_u)^2,\\
\sigma^2(u(t-s))=\sigma^2(ut_u)+\dot{\sigma^2}(ut_u)u(t-t_u-s)+\frac{1}{2}\ddot{\sigma^2}(u\theta_{2,u})u^2(t-t_u-s)^2,
\EQNY
with $\theta_{1,u}\in (t,t_u)$ and $\theta_{2,u}\in (t-s,t_u)$.
Inserting the above expansions to the following equation  and using (\ref{eq4}), we have
\BQN\label{eq6}
1-\frac{\sigma_{\gamma,u}^2(s,t)}{\sigma_{\gamma,u}^2(0,t_u)}&=&\frac{(1+c(t-\gamma s))^2\sigma^2(ut_u)-\left[(1-\gamma)\sigma^2(ut)+(\gamma^2-\gamma)\sigma^2(us)+\gamma\sigma^2(u(t-s))\right](1+ct_u )^2}{(1+c(t-\gamma s))^2\sigma^2(ut_u)}\nonumber\\
&=&\frac{(\gamma-\gamma^2)(1+ct_u)^2\sigma^2(us)-\frac{1-\gamma}{2}u^2\ddot{\sigma^2}(u\theta_{1,u})(1+ct_u)^2(t-t_u)^2}{(1+c(t-\gamma s))^2\sigma^2(ut_u)}\nonumber\\
&&+\frac{\sigma^2(ut_u)c^2(t-t_u-\gamma s)^2-\frac{\gamma}{2}\ddot{\sigma^2}(u\theta_{2,u})u^2(1+ct_u)^2(t-t_u-s)^2}{(1+c(t-\gamma s))^2\sigma^2(ut_u)} \nonumber\\
&=&\frac{(\gamma-\gamma^2)(1+ct_u)^2\sigma^2(us)+\left(\sigma^2(ut_u)c^2-\frac{1-\gamma}{2}u^2\ddot{\sigma^2}(u\theta_{1,u})(1+ct_u)^2-\frac{\gamma}{2}\ddot{\sigma^2}(u\theta_{2,u})u^2(1+ct_u)^2\right)(t-t_u)^2}{(1+c(t-\gamma s))^2\sigma^2(ut_u)}\nonumber\\
&&+\frac{\sigma^2(ut_u)c^2(\gamma^2 s^2-2\gamma(t-t_u)s)-\frac{\gamma}{2}\ddot{\sigma^2}(u\theta_{2,u})u^2(1+ct_u)^2(s^2-2(t-t_u)s)}{(1+c(t-\gamma s))^2\sigma^2(ut_u)}
\EQN
 \rd{It follows from (\ref{re3}) that }for  any $\delta>0$ and $u$ large enough
\BQN\label{eq5}
\frac{s^2}{\frac{\sigma^2(us)}{\sigma^2(ut_u)}}=t_u^2\frac{g_2(ut_u)}{g_2(us)}\leq 2t_*^{2\alpha_\IF} \delta^{2-2\alpha_\IF},\ \ s\in (0,\delta].
\EQN
Following (\ref{re2}), we have that
 \BQN\label{eq7}
 \frac{\ddot{\sigma^2}(u\theta_{i,u})u^2}{\sigma^2(ut_u)}\sim \frac{2\alpha_\IF(2\alpha_\IF-1)}{(t_*)^2}, \quad u\rw\IF, i=1,2.
 \EQN
 Moreover, for $\delta>0$ sufficiently small and $u$ sufficiently large
 \BQN\label{eq8}
 |t-t_u| s&\leq&  \delta^{(1-\alpha_\IF)/2}|t-t_u| \delta^{-(1-\alpha_\IF)/2}s\nonumber\\
 &\leq& \delta^{1-\alpha_\IF}(t-t_u)^2+ \delta^{\alpha_\IF-1}s^2\notag \\
 &\leq &  Q \delta^{1-\alpha_\IF}\left(\frac{\sigma^2(us)}{\sigma^2(u)}+(t-t_u)^2\right), \quad s\in (0,\delta].
 \EQN
 Hence inserting (\ref{eq5})-(\ref{eq8}) into (\ref{eq6}), we have that for $u$ sufficiently large
\BQNY
 (2a_1-Q \delta^{1-\alpha_\IF})(t-t_u)^2+(2a_2-Q \delta^{1-\alpha_\IF})\frac{\vf (us)}{\vf (u)}&\leq& 1-\frac{\sigma_{\gamma,u}^2(s,t)}{\sigma_{\gamma,u}^2(0,t_u)}\\
& \leq& (2a_1+Q \delta^{1-\alpha_\IF})(t-t_u)^2+(2a_2+Q \delta^{1-\alpha_\IF})\frac{\vf (us)}{\vf (u)}
\EQNY
 for $0<s< \delta$ and $|t-t_u|<\delta$ with $\delta>0$ sufficiently small and  $Q$ a fixed constant, which establishes the claim.
\QED
\\
\prooflem{L2}. It follows from the direct calculation that
\BQNY
1-r_u(s,t,s_1,t_1)=\frac{D_{1,u}(s,t,s_1,t_1)-D_{2,u}(s,t,s_1,t_1)+\gamma D_{3,u}(s,t,s_1,t_1)}{2\sigma_\gamma(us,ut)\sigma_\gamma(us_1,ut_1)},
\EQNY
with
\BQNY
D_{1,u}(s,t,s_1,t_1)&=&\sigma^2(u|t-t_1|)+\gamma^2\sigma^2(u|s-s_1|), \ \ D_{2,u}(s,t,s_1,t_1)=(\sigma_\gamma(us,ut)-\sigma_\gamma(us_1,ut_1))^2,\\
D_{3,u}(s,t,s_1,t_1)&=&\sigma^2(u|t-s|)+\sigma^2(u|t_1-s_1|)-\sigma^2(u|t_1-s|)-\sigma^2(u|t-s_1|).
\EQNY
Using Taylor expansion, we have
\BQNY
D_{3,u}(s,t,s_1,t_1)
&=&u\dot{\sigma^2}(u(t_1-s))(t-t_1)+\frac{1}{2}u^2\ddot{\sigma^2}(u\theta_1)(t-t_1)^2\\
&&+u\dot{\sigma^2}(u(t-s_1))(t_1-t)
+\frac{1}{2}u^2\ddot{\sigma^2}(u\theta_2)(t-t_1)^2\\
&=&\frac{1}{2}u^2\ddot{\sigma^2}(u\theta_1)(t-t_1)^2+\frac{1}{2}u^2\ddot{\sigma^2}(u\theta_2)(t-t_1)^2\\
&&+u^2\ddot{\sigma^2}(u\theta_3)(t-t_1)(t_1-t+s_1-s)\\
&\leq&u^2\left(\frac{1}{2}\ddot{\sigma^2}(u\theta_1)+\frac{1}{2}\ddot{\sigma^2}(u\theta_2)+2\ddot{\sigma^2}(u\theta_3)\right)(t-t_1)^2+2u^2\ddot{\sigma^2}(u\theta_3)(s-s_1)^2,
\EQNY
where $\theta_1$, $\theta_2$ and $\theta_3$ are some positive constants (depending on $u$) satisfying $\frac{t_*}{2}<\theta_i<\frac{3}{2}t_*, i=1,2,3$  for $u$ sufficiently large.
\rd{From (\ref{re2}) and (\ref{A1})},  we have that for $\delta>0$
\BQNY
\sup_{t\in(0,\delta)}\Big\vert\frac{u^2\ddot{\sigma^2}(u)t^2}{\sigma^2(ut)}\Big\vert\leq Q\sup_{t\in(0,\delta)}\frac{\sigma^2(u)t^2}{\sigma^2(ut)}=Q\sup_{t\in(0,\delta)}\frac{g_2(u)}{g_2(ut)},
\EQNY
which together with (\ref{re3}) implies that if $\delta_u\rw 0$ as $u\rw\IF$
\BQNY
\sup_{t\in(0,\delta_u)}\Big\vert\frac{u^2\ddot{\sigma^2}(u)t^2}{\sigma^2(ut)}\Big\vert\leq Q\sup_{t\in(0,\delta_u)}\frac{g_2(u)}{g_2(ut)}\rw 0, \ \ u\rw\IF.
\EQNY
Therefore we get that  uniformly for $(s,t)\neq (s_1,t_1)\in [0,\delta_u)\times(t_u-\delta_u,t_u+\delta_u)$
\BQNY
\frac{D_{3,u}(s,t,s_1,t_1)}{D_{1,u}(s,t,s_1,t_1)}\rw 0, \ \ u\rw\IF.
\EQNY
By (\ref{A1}) and {\bf AIII} we have for any $x\in(0,\IF)$ and any $y\in[0,1]$
\BQNY
1\geq \frac{\sigma^2(xy)}{\sigma^2(x)}=\frac{g_2(xy)}{g_2(x)}y^2\geq y^2.
\EQNY
Hence by  UCT for $0\leq s_1<s<\delta_u$ with $\delta_u\rw 0$
\BQN\label{eq21}
\frac{\left(\sigma^2(us)-\sigma^2(us_1)\right)^2}{\sigma^2(u|s-s_1|)\sigma^2(u)}=\frac{\sigma^2(us)}{\sigma^2(u)}
\frac{\left(1-\frac{\sigma^2(us_1)}{\sigma^2(us)}\right)^2}{\frac{\sigma^2(us(1-s_1/s))}{\sigma^2(us)}}\leq \frac{\sigma^2(us)}{\sigma^2(u)}(1+s_1/s)^2\leq 4\frac{\sigma^2(us)}{\sigma^2(u)}\rw 0, \ \ u\rw\IF.
\EQN
By (\ref{re1}) and (\ref{A1}) we have
\BQNY
\lefteqn{\frac{D_{2,u}(s,t,s_1,t_1)}{D_{1,u}(s,t,s_1,t_1)}}\\
&\leq & 4\frac{(1-\gamma)^2(\sigma^2(ut)-\sigma^2(ut_1))^2+\gamma^2(\sigma^2(u(t-s))-\sigma^2(u(t_1-s_1)))^2+(\gamma-\gamma^2)^2(\sigma^2(us)-\sigma^2(us_1))^2}
{D_{1,u}(s,t,s_1,t_1)\left(\sigma_\gamma(us,ut)+\sigma_\gamma(us_1,ut_1)\right)^2}\\
&\leq &Q\left(\frac{(u\dot{\sigma^2}(u))^2(t-t_1)^2}{\sigma^2(u)\sigma^2(u|t-t_1|)}+
\frac{(u\dot{\sigma^2}(u))^2(s-s_1)^2}{\sigma^2(u)\sigma^2(u|s-s_1|)}+\frac{\left(\sigma^2(us)-\sigma^2(us_1)\right)^2}{\sigma^2(u|s-s_1|)\sigma^2(u)}\right)\\
&\leq & Q_1\left(\frac{g_2(u)}{g_2(u|t-t_1|)}+\frac{g_2(u)}{g_2(u|s-s_1|)}+\frac{\left(\sigma^2(us)-\sigma^2(us_1)\right)^2}
{\sigma^2(u|s-s_1|)\sigma^2(u)}\right).
\EQNY
Further, it follows from (\ref{re3}) and (\ref{eq21}) that
\BQNY
\frac{D_{2,u}(s,t,s_1,t_1)}{D_{1,u}(s,t,s_1,t_1)}
\rw 0,
\EQNY
as $u\rw\IF$ uniformly for $(s,t)\neq (s_1,t_1)\in [0,\delta_u)\times(t_u-\delta_u,t_u+\delta_u)$ with $\delta_u\rw 0$.  This completes the proof. \QED
\\
\prooflem{Var} \CE{We have}
\BQNY
\sigma_{1,u}^2(s,t)&=&\frac{(1-\gamma)\sigma^2(t)+(\gamma^2-\gamma)\sigma^2(s)+\gamma\sigma^2(t-s)}{\sigma^2(T)}\frac{(u+cT)^2}{(u+c(t-\gamma s))^2}\\
&=:&f_1(s,t) f_{2,u}(s,t),\ \ (s,t)\in D_T=\{(s,t), 0\leq s\leq t\leq T\}.
\EQNY
In light of {\bf BIII}, $f_1(s,t)$ is strictly increasing with respect to $t$ and strictly decreasing with respect to $s$ for $(s,t)\in D_T$. Moreover,
$$\lim_{u\rw\IF}\sup_{(s,t)\in D_T}|f_{2,u}(s,t)-1|=0.$$
Thus we conclude that the maximum value of $\sigma_{1,u}^2(s,t)$ over $D_T$ must be attained in a sufficiently small neighbourhood of $(0,T)$ for $u$ large enough. Further, as $(s,t)\rw (0,T)$
 $$ 1-f_1(s,t)=\frac{\vfp(T)}{\vf(T)}(T-t)(1+o(1))+\left\{
\begin{array}{cc}
\frac{\gamma\vfp(T)}{\vf(T)}s(1+o(1)),& \text{   if  }\vf(s)=o(s),\\
\frac{b(\gamma-\gamma^2)+\gamma\vfp(T)}{\vf(T)}s(1+o(1)),& \text{  if  }\vf(s)\sim b s,\\
\frac{\gamma-\gamma^2}{\vf(T)}\vf(s)(1+o(1)),& \text{  if  }s=o(\vf(s)),
 \end{array}
 \right.$$
 and for $u>1$
 $$1-f_{2,u}(s,t)=\frac{-2c}{u+cT}\left(T-t+\gamma s\right)(1+o(1)),$$
 which imply that (\ref{Vare}) holds and further the maximum point of $\sigma_{1,u}(s,t)$ in a neighbourhood of $(0,T)$ is $(0,T)$. Thus the claim is established. \QED

\prooflem{CorF} The proof is similar to that of Lemma \ref{L2}. We have
\BQNY
1-r_1(s,t,s_1,t_1)=\frac{D_{1}(s,t,s_1,t_1)-D_{2}(s,t,s_1,t_1)+\gamma D_{3}(s,t,s_1,t_1)}{2\sigma_\gamma(s,t)\sigma_\gamma(s_1,t_1)},
\EQNY
with
\BQNY
D_{1}(s,t,s_1,t_1)&=&\sigma^2(|t-t_1|)+\gamma^2\sigma^2(|s-s_1|), \ \ D_{2}(s,t,s_1,t_1)=(\sigma_\gamma(s,t)-\sigma_\gamma(s_1,t_1))^2,\\
D_{3}(s,t,s_1,t_1)&=&\sigma^2(|t-s|)+\sigma^2(|t_1-s_1|)-\sigma^2(|t_1-s|)-\sigma^2(|t-s_1|).
\EQNY
Using Taylor expansion and the fact that $t^2=o(\sigma^2(t))$ as $t\downarrow 0$, we have
\BQNY
D_{3}(s,t,s_1,t_1)
&=&\dot{\sigma^2}(t_1-s)(t-t_1)+\frac{1}{2}\ddot{\sigma^2}(\theta_4)(t-t_1)^2
+\dot{\sigma^2}(t-s_1)(t_1-t)
+\frac{1}{2}\ddot{\sigma^2}(\theta_5)(t-t_1)^2\\
&=&\frac{1}{2}\ddot{\sigma^2}(\theta_4)(t-t_1)^2+\frac{1}{2}\ddot{\sigma^2}(\theta_5)(t-t_1)^2
+\ddot{\sigma^2}(\theta_6)(t-t_1)(t_1-t+s_1-s)\\
&\leq&\left(\frac{1}{2}\ddot{\sigma^2}(\theta_4)+\frac{1}{2}\ddot{\sigma^2}(\theta_5)+2\ddot{\sigma^2}(\theta_6)\right)(t-t_1)^2
+2\ddot{\sigma^2}(\theta_6)(s-s_1)^2\\
&=& o\left(D_{1}(s,t,s_1,t_1)\right), \ \ s,s_1\rw 0, t,t_1\rw T,
\EQNY
where $\theta_4$, $\theta_5$ and $\theta_6$ are some positive constants satisfying $\frac{T}{2}<\theta_i<\frac{3}{2}T, i=4,5,6$.
By (\ref{A1}) and {\bf BIII} we have for any $x\in(0,\IF)$ and any $y\in[0,1]$
\BQNY
1\geq \frac{\sigma^2(xy)}{\sigma^2(x)}=\frac{g_2(xy)}{g_2(x)}y^2\geq y^2,
\EQNY
hence for $0\leq s_1<s<T/2$
\BQN\label{eq211}
\frac{\left(\sigma^2(s)-\sigma^2(s_1)\right)^2}{\sigma^2(|s-s_1|)}&=& \sigma^2(s)
\frac{\left(1-\frac{\sigma^2(s_1)}{\sigma^2(s)}\right)^2}{\frac{\sigma^2(s(1-s_1/s))}{\sigma^2(s)}} \notag \\
&\leq & \sigma^2(s)(1+s_1/s)^2\notag \\
&\leq & 4\sigma^2(s)\rw 0, \ \ s\rw 0.
\EQN
By (\ref{A1}), (\ref{eq211}) and the fact that $t^2=o(\sigma^2(t))$ as $t\downarrow 0$, we have
\BQNY
D_{2}(s,t,s_1,t_1)&=&\frac{(\sigma_\gamma^2(s,t)-\sigma_\gamma^2(s_1,t_1))^2}{(\sigma_\gamma(s,t)+\sigma_\gamma(s_1,t_1))^2}\\
&=&\frac{\left((1-\gamma)(\sigma^2(t)-\sigma^2(t_1))+(\gamma^2-\gamma)(\sigma^2(s)-\sigma^2(s_1))+\gamma(\sigma^2(t-s)-\sigma^2(t_1-s_1
))\right)^2}{(\sigma_\gamma(s,t)+\sigma_\gamma(s_1,t_1))^2}\\
&\leq& \frac 8 {\sigma^2(T)} \left((\dot{\sigma^2}(T))^2(t-t_1)^2+(\dot{\sigma^2}(T))^2(t-t_1-s+s_1)^2+(\sigma^2(s)-\sigma^2(s_1))^2\right)\\
&=& o\left(D_{1}(s,t,s_1,t_1)\right), \ \ s,s_1\rw 0, t,t_1\rw T.
\EQNY
Therefore, we have
$$1-r_1(s,t,s_1,t_1)\sim \frac{\sigma^2(|t-t_1|)+\gamma^2\sigma^2(|s-s_1|)}{2\sigma^2(T)}, \ \ s,s_1\rw 0, t,t_1\rw T,$$
 which completes the proof. \QED

{\bf Acknowledgement}: \eG{We \kk{thank referees} for several suggestions and comments.} Thanks to Swiss National Science Foundation grant No.200021-166274.
KD acknowledges partial support by NCN Grant No 2015/17/B/ST1/01102 (2016-2019).

 \bibliographystyle{ieeetr}

 \bibliography{gamE}

\end{document}